\documentclass[11pt]{article}
\usepackage[a4paper,margin=1in]{geometry}
\usepackage{amsmath,amssymb,amsthm,mathtools,bm,graphicx}
\usepackage{booktabs,enumitem,microtype,aliascnt,array}
\usepackage[colorlinks=true,linkcolor=blue,citecolor=blue,urlcolor=blue]{hyperref}
\usepackage[nameinlink,noabbrev]{cleveref}
\hypersetup{pdftitle={Source Elimination, Covariance Geometry, and Calibration for Inverse Problems with Latent Perturbations},pdfauthor={Saman Khoramian}}
\numberwithin{equation}{section}
\newtheorem{assumption}{Assumption}[section]
\newaliascnt{theorem}{assumption}
\newtheorem{theorem}[theorem]{Theorem}
\aliascntresetthe{theorem}
\newaliascnt{proposition}{assumption}
\newtheorem{proposition}[proposition]{Proposition}
\aliascntresetthe{proposition}
\newaliascnt{corollary}{assumption}
\newtheorem{corollary}[corollary]{Corollary}
\aliascntresetthe{corollary}
\newaliascnt{lemma}{assumption}
\newtheorem{lemma}[lemma]{Lemma}
\aliascntresetthe{lemma}
\newaliascnt{remark}{assumption}
\newtheorem{remark}[remark]{Remark}
\aliascntresetthe{remark}
\crefname{assumption}{assumption}{assumptions}
\Crefname{assumption}{Assumption}{Assumptions}
\crefname{theorem}{theorem}{theorems}
\Crefname{theorem}{Theorem}{Theorems}
\crefname{proposition}{proposition}{propositions}
\Crefname{proposition}{Proposition}{Propositions}
\crefname{corollary}{corollary}{corollaries}
\Crefname{corollary}{Corollary}{Corollaries}
\crefname{lemma}{lemma}{lemmas}
\Crefname{lemma}{Lemma}{Lemmas}
\crefname{remark}{remark}{remarks}
\Crefname{remark}{Remark}{Remarks}
\newcommand{\R}{\mathbb R}
\newcommand{\cA}{\mathcal A}
\newcommand{\cB}{\mathcal B}
\newcommand{\cC}{\mathcal C}
\newcommand{\cE}{\mathcal E}

\newcommand{\cK}{\mathcal K}

\newcommand{\cQ}{\mathcal Q}
\newcommand{\cR}{\mathcal R}
\newcommand{\cU}{\mathcal U}
\newcommand{\ip}[2]{\langle #1,#2\rangle}
\newcommand{\norm}[1]{\lVert #1\rVert}
\newcommand{\abs}[1]{\lvert #1\rvert}
\newcommand{\tr}{\operatorname{tr}}
\newcommand{\argmin}{\operatorname*{arg\,min}}
\newcommand{\PSD}{\mathbb S_+}

\newcommand{\E}{\mathbb E}
\newcommand{\Cov}{\operatorname{Cov}}
\newcommand{\Var}{\operatorname{Var}}
\newcommand{\Ran}{\operatorname{ran}}
\newcommand{\Ker}{\operatorname{ker}}
\title{Source Elimination, Covariance Geometry, and Calibration for Inverse Problems with Latent Perturbations}
\author{Saman Khoramian\\Kharazmi University}
\date{November 24, 2026}
\begin{document}
\maketitle

\begin{abstract} 
We study the interface between quadratic source relaxation, measurement-space geometry, and independent replicate calibration in PDE inverse problems. Classical Hilbert-space elimination reduces a coercively penalized latent source to a finite-dimensional quadratic involving $\Sigma+T_cC_{\mathrm{pen},\theta}(c)T_c^*$. We distinguish this deterministic profile from Gaussian marginal likelihood, constrained source fitting, and local statistical optimality.

The calibration target is a structured compression $K=G+GSG$. For IID symmetric matrix-valued replicates with mean $G$ and finite second moments, an order-two cross-replicate estimator is unbiased and admits an explicit Hoeffding variance decomposition. For fixed \(S\), or conditionally on an independently learned \(S\), positive-semidefinite projection gives an admissible inverse weight with a mean-square calibration bound of order \(n^{-1}\); unconditional bounds for random \(S\) require suitable integrability. The same rate transfers to the inverse weight; under an additional strong-convexity assumption, it also controls reconstruction error relative to an oracle frozen-weight estimator, not relative to the true coefficient.

A one-dimensional reciprocal elliptic example provides an explicit excitation--receiver protocol for the matrix replicates. Its Green matrix is an observation-map energy Gram, not the Gram of an arbitrary forcing-to-data map. A fixed-reference Hilbert formulation makes this distinction and its coefficient derivatives precise. Candidate-dependent exact relaxation and reference-frozen calibration are defined as different estimators, with matched training objectives. Compactness-based optimistic bilevel existence and local smooth-branch hypergradients are stated separately.

This is a theory note with three limited numerical illustrations: a synthetic elliptic calibration experiment, a hypothetical seismic-motivated RMSE-margin calculation, and a synthetic one-dimensional Helmholtz calibration experiment. They are not empirical reconstruction studies and provide no evidence of practical superiority. It establishes neither physical covariance recovery nor general improvements in nonlinear inverse-problem optimization.
\end{abstract}
\noindent\textbf{Keywords:} inverse problems; source relaxation; covariance geometry; bilevel learning; calibration; U-statistics.

\section{Introduction}\label{sec:intro}
\subsection{The reduced object}\label{sec:intro-reduced}
An uncertain excitation or model discrepancy acts before observation, whereas measurement noise acts in data space. A quadratic relaxation distinguishes these roles. For a bounded map $T_c:\cQ\to\R^m$, a positive-definite measurement weight $\Sigma$, and a coercive source geometry $C$, the problem
\[
 \min_{q\in\cQ}\left\{\frac12\norm{T_cq-e(c)}_{\Sigma^{-1}}^2
                +\frac12\ip{C^{-1}q}{q}_{\cQ}\right\}
\]
has minimizer and value
\[
 q^\star=CT_c^*(\Sigma+T_cCT_c^*)^{-1}e(c),\qquad
 \Phi(c)=\frac12 e(c)^T(\Sigma+T_cCT_c^*)^{-1}e(c).
\]
Here $\norm{x}_{W}^2:=x^TWx$ for a positive-definite matrix $W$, with the analogous Hilbert quadratic-form notation for operators. The reduced objective sees the source geometry only through its compression $K=T_cCT_c^*$. It does not, in general, identify $C$, and it does not turn the resulting coefficient objective into a convex function.

This note makes precise a particular learning--calibration interface. A supervised stage selects a geometry parameter. Independent matrix-valued replicates then estimate a structured compression at one physical configuration. Using that fixed estimate throughout coefficient optimization produces a \emph{frozen-weight estimator}; it is not the exact candidate-dependent source profile above. Keeping these estimators distinct is essential to the interpretation of the results.

\subsection{Classical ingredients and the specific increment}\label{sec:intro-related}
Quadratic elimination is classical Schur/Woodbury algebra, equivalently the dual form of quadratic Tikhonov regularization \cite{GolubVanLoan2013,EnglHankeNeubauer1996}. With
\[
 A=\Sigma+TCT^*,\qquad H=C^{-1}+T^*\Sigma^{-1}T,
\]
the correctly weighted identities are
\begin{equation}\label{eq:push-through}
 A^{-1}=\Sigma^{-1}-\Sigma^{-1}TH^{-1}T^*\Sigma^{-1},\qquad
 A^{-1}TC=\Sigma^{-1}TH^{-1}.
\end{equation}
One can solve in data space or source space, depending on dimension, available operator applications, and conditioning. No computational advantage follows from the identity alone.

\paragraph{Correspondence with state relaxation.}
For an invertible PDE operator $\cA(c):\cU\to\cU^*$, fixed forcing $f\in\cU^*$, and $\cQ=\cU^*$ equipped with a specified fixed Hilbert structure, consider
\[
 \min_{u,q}\left\{\tfrac12\ip{C^{-1}q}{q}_{\cQ}
              +\tfrac12\norm{Mu-d}_{\Sigma^{-1}}^2\right\},
 \qquad \cA(c)u=f+q.
\]
Substitution of $u=\cA(c)^{-1}(f+q)$ gives exactly the source problem with $T_c=M\cA(c)^{-1}$ and $e(c)=d-M\cA(c)^{-1}f$. Thus this forcing relaxation is a reparameterization of a state-relaxation problem, not a different elimination principle \cite{vanLeeuwenHerrmann2016,vanLeeuwenYang2025}. The reciprocal example below uses a separately specified state-perturbation geometry; it does not infer its Green identity from this forcing model.

\paragraph{Closest PDE antecedent.}
Van Leeuwen and Yang \cite{vanLeeuwenYang2025} already reduce a relaxed PDE-constrained problem to weighted nonlinear least squares, with an inverse weight built from a Gram matrix of PDE solutions. They also show how to derive data-driven Grams in several settings. In their one-dimensional elliptic example (Section~3.1.2), colocated point sources and receivers give the same response-matrix/energy-Gram identity used in \cref{lem:reciprocity}. Neither this identity nor the weighted-residual reduction is claimed as new here. We use the same idealized forward responses, but stipulate additional, independently repeated full-matrix sessions at a fixed physical coefficient, with centered errors and explicit acquisition cost (the protocol in \cref{sec:stability-reciprocal}). This differs from constructing a Gram from an available measurement set: the repeated sessions supply IID noisy matrix replicates rather than requiring the ordinary reconstruction data to supply the Gram.

The specific increment is the calibration of the structured, generally noncommuting target $G+GSG$ from these replicates. The distinct-replicate estimator is unbiased, with an exact second-moment variance decomposition. The resulting mean-square bound transfers through PSD projection to an admissible inverse weight and, under additional strong convexity, to the distance from an oracle frozen coefficient estimator (\cref{thm:calibration,cor:calibration-consistency,prop:frozen-stability}).

Bilevel existence and sensitivity are applications of parametric optimization \cite{BonnansShapiro2000,DontchevRockafellar2009,Dempe2002}, and supervised variational learning has substantial precedent \cite{DeLosReyesSchonlieb2013,DeLosReyesSchonliebValkonen2016,DeLosReyesSchonliebValkonen2017,GuntherKroner2021,FranceschiEtAl2018}. The earlier preprint \cite{Khoramian2020} provides a related supervised regularization template, not a novelty argument for the present algebra. Gaussian conditioning and covariance propagation are also classical \cite{RasmussenWilliams2006,Stuart2010,KaipioSomersalo2005}.

The calibration estimator is a matrix-valued specialization of order-two U-statistics \cite{Hoeffding1948,Serfling1980}. Its useful application-specific content is the pair of noncommutative kernels associated with $G+GSG$, the explicit finite-sample constants under only second moments, and their combination with an admissibility rule and a specified reconstruction estimator. We do not claim a new general U-statistic theorem, a new elimination identity, or a new generic hypergradient theorem.

\begin{table}[t]
\centering\small
\caption{Classical tools and their role in this note.}
\label{tab:positioning}
\begin{tabular}{@{}>{\raggedright\arraybackslash}p{0.23\textwidth}>{\raggedright\arraybackslash}p{0.31\textwidth}>{\raggedright\arraybackslash}p{0.36\textwidth}@{}}
\toprule
Tool & Established ingredient & Specialization/interface here \\
\midrule
Quadratic elimination & Schur/Woodbury identity & Coercive versus factor-domain geometry; candidate-dependent versus frozen estimator \\
Order-two U-statistics & Hoeffding orthogonality & Explicit kernels for $G+GSG$; PSD and inverse-weight error bounds \\
Parametric optimization & Compact-graph existence and local KKT sensitivity & Matched lower objectives; full geometry derivatives; conditional reconstruction bound \\
\bottomrule
\end{tabular}
\end{table}

\subsection{Results and scope}\label{sec:intro-contrib}
The main statements are the precise deterministic reduction and its semidefinite extension (\cref{thm:representer,prop:factor}), an acquisition-compatible reciprocal example (\cref{lem:reciprocity}), and the calibration/inverse-weight/reconstruction bounds (\cref{thm:calibration,cor:calibration-consistency,prop:frozen-stability}). The last reconstruction bound assumes strong convexity of an oracle frozen objective. It does not derive that property from norm learning.

The numerical material comprises a synthetic elliptic calibration experiment, a hypothetical seismic-motivated RMSE-margin illustration, and a synthetic one-dimensional Helmholtz calibration experiment. These examples test distinct calibration and stability mechanisms under controlled assumptions; they are not empirical PDE reconstruction studies and do not establish practical superiority. Supplementary code checks algebra and finite-support expectations (\cref{app:verification}). We claim neither population consistency of supervised learning nor global nonlinear identifiability, physical covariance recovery, or transfer of a calibration matrix to another physical coefficient. Calibration consistency for a fixed target is proved and must not be confused with these stronger claims.

\section{Mathematical setting and notation}\label{sec:setting}
\subsection{Fixed Hilbert structures and coefficient-dependent geometries}\label{sec:setting-spaces}
All deterministic Hilbert spaces are real. Let $\cQ$ be a fixed Hilbert space and $\cC$ a coefficient set in a normed space. For each $c\in\cC$, let $T_c:\cQ\to\R^m$ be bounded and let
\begin{equation}\label{eq:residual}
 e_j(c):=d_j-\mu_j(c)\in\R^m
\end{equation}
be a residual relative to a modeled mean. Hilbert adjoints are always taken in the declared \emph{fixed} inner products unless an energy adjoint is explicitly marked.

A standard forcing specialization is
\begin{equation}\label{eq:T}
 T_c=M\cA(c)^{-1}\cB,
\end{equation}
where $\cU$ is fixed, $\cA(c):\cU\to\cU^*$ is boundedly invertible, $\cB:\cQ\to\cU^*$ and $M:\cU\to\R^m$ are fixed bounded operators, and $\mu(c)=M\cA(c)^{-1}f$ for fixed $f$. A second specialization, defined in \cref{sec:stability-reciprocal}, uses a state perturbation, $T=M$, and places the coefficient dependence in the source geometry. This is not obtained from \cref{eq:T} by treating an unspecified injection as a Riesz map.

\begin{assumption}[Operator regularity]\label{ass:operator}
The maps $c\mapsto T_c$ and $c\mapsto e_j(c)$ are norm-continuous. When derivatives are used, the relevant maps extend continuously Fr\'echet differentiably to an open neighborhood. For \cref{eq:T}, these properties follow from the corresponding regularity of $\cA$, continuity of the modeled means, and a uniform bound on $\norm{\cA(c)^{-1}}$ on the set under consideration.
\end{assumption}
Continuity alone does not imply local Lipschitz continuity. The continuously differentiable case does imply it locally. Since the data space is finite-dimensional, every bounded $T_c$ has finite rank and is compact, but it need not be injective.

\begin{assumption}[Coercive deterministic geometry]\label{ass:geometry}
The operators $C_{\mathrm{pen},\theta}(c):\cQ\to\cQ$ are bounded and self-adjoint, and
\[
 0<\kappa_- I\preceq C_{\mathrm{pen},\theta}(c)\preceq\kappa_+ I<\infty
\]
uniformly over the pairs $(c,\theta)$ under consideration. The measurement weight $\Sigma=\Sigma_\theta$ is positive definite and is independent of the candidate coefficient $c$.
\end{assumption}
The bounds are convenient sufficient assumptions for a bounded precision and a strongly convex source problem, not necessary conditions for every elimination formula. Semidefinite geometries require the range-domain formulation in \cref{prop:factor}.

\subsection{Penalty geometry is not automatically a covariance}\label{sec:setting-geometry}
We reserve $C_{\mathrm{pen},\theta}$ for deterministic geometry and $C_{\mathrm{prob}}$ for a probabilistic source covariance. On a separable infinite-dimensional Hilbert space, a Gaussian covariance is positive and trace class. A coercive bounded operator is not trace class (\cref{prop:trace-class}); thus these objects cannot simply be identified there.

For a positive trace-class covariance $C$, its Cameron--Martin space is $\Ran(C^{1/2})$ with norm
\[
 \norm{h}_{C}:=\min\{\norm{z}_{\cQ}:C^{1/2}z=h\}.
\]
The minimizing preimage is in $(\Ker C)^\perp$. If $C$ is injective, this is the usual $\norm{C^{-1/2}h}_{\cQ}$. It is not the ambient Hilbert norm, and typical Gaussian draws need not belong to this space.

\begin{table}[t]
\centering\small
\caption{Main notation. Operator adjoints refer to fixed Hilbert structures.}
\label{tab:notation}
\begin{tabular}{@{}ll@{}}
\toprule
Symbol & Meaning \\
\midrule
$m,p,d_\theta$ & data dimension, discretized coefficient dimension, parameter dimension \\
$n,n_{\mathrm{tr}}$ & calibration replicate count, training sample count \\
$T_c,e(c)$ & latent-perturbation map and data residual \\
$C_{\mathrm{pen},\theta}(c),C_{\mathrm{prob}}$ & deterministic geometry, probabilistic covariance \\
$G(c)=T_cT_c^*$ & fixed-space identity-geometry Gram \\
$C_0(c),G_0(c)=D(c)$ & reciprocal model's base geometry and energy Gram \\
$K_\theta(c),S_\theta$ & source compression and its structured matrix parameter \\
$\Sigma,\varrho,\cR$ & measurement weight, regularization weight, coefficient penalty \\
$\widehat K_U,\widehat K_+$ & raw cross-replicate estimate and PSD projection \\
$\Var_F(Z\mid S)$ & conditional squared-Frobenius variance \\
$J,P,\Gamma^\dagger$ & local Jacobian, coefficient covariance, true local noise covariance \\
$\mu(c)$ & modeled data mean \\
$\mu_{\mathrm{sc}}$ & strong-convexity modulus \\
\bottomrule
\end{tabular}
\end{table}

\section{Exact source elimination}\label{sec:elimination}
\subsection{Coercive geometry and the profile objective}\label{sec:elimination-representer}
Fix $c$, $e=e(c)$, and $\theta$. Suppress the arguments of $C=C_{\mathrm{pen},\theta}(c)$ and $T=T_c$ when convenient. Define
\begin{equation}\label{eq:source-objective}
 \cE_\theta(q;c)=\frac12\norm{T_cq-e(c)}_{\Sigma^{-1}}^2
                  +\frac12\ip{C_{\mathrm{pen},\theta}(c)^{-1}q}{q}_{\cQ}.
\end{equation}
\begin{theorem}[Exact source elimination]\label{thm:representer}
Let $T_c$ be bounded and let \cref{ass:geometry} hold. The source objective is coercive and strongly convex with modulus $\kappa_+^{-1}$. Its unique minimizer is
\begin{equation}\label{eq:qstar}
 q_\theta^\star(c)=C_{\mathrm{pen},\theta}(c)T_c^*
                 (\Sigma+K_\theta(c))^{-1}e(c),\qquad
 K_\theta(c)=T_cC_{\mathrm{pen},\theta}(c)T_c^*,
\end{equation}
and its minimum is
\begin{equation}\label{eq:profile}
 \Phi_\theta(c)=\frac12e(c)^T(\Sigma+K_\theta(c))^{-1}e(c).
\end{equation}
Writing $r=e-T_cq_\theta^\star=\Sigma(\Sigma+K_\theta)^{-1}e$, the energy split is
\begin{equation}\label{eq:energy-split}
 \Phi_\theta(c)=\tfrac12\norm{r}_{\Sigma^{-1}}^2
               +\tfrac12\norm{q_\theta^\star}_{C_{\mathrm{pen},\theta}(c)^{-1}}^2.
\end{equation}
\end{theorem}
\begin{proof}
The Hessian is $H=C^{-1}+T^*\Sigma^{-1}T\succeq\kappa_+^{-1}I$. Hence the objective has a unique minimizer satisfying $Hq=T^*\Sigma^{-1}e$. Put $A=\Sigma+TCT^*$ and $q=CT^*A^{-1}e$. Since $Tq=e-\Sigma A^{-1}e$, direct substitution gives
\[
 T^*\Sigma^{-1}(Tq-e)+C^{-1}q=-T^*A^{-1}e+T^*A^{-1}e=0.
\]
This proves \cref{eq:qstar} without any unweighted push-through identity. At this source the two energy terms are $\tfrac12e^TA^{-1}\Sigma A^{-1}e$ and $\tfrac12e^TA^{-1}TCT^*A^{-1}e$. Their sum is $\tfrac12e^TA^{-1}e$.
\end{proof}

\begin{remark}[Compression versus the source]\label{rem:elimination-scope}
The profile depends on $C$ only through $TCT^*$, while the minimizing source also depends on $CT^*$. Equal compressions give equal profile objectives but need not give equal sources. This is a deterministic statement, with no covariance interpretation required.
\end{remark}

\begin{corollary}[Coefficient profile and attainment]\label{cor:profile}
Let $\cR:\cC\to\R$ and $\varrho\ge0$. Pointwise source elimination gives
\begin{equation}\label{eq:reduced}
 \inf_{c\in\cC}\min_{q\in\cQ}\{\cE_\theta(q;c)+\varrho\cR(c)\}
 =\inf_{c\in\cC}\{\Phi_\theta(c)+\varrho\cR(c)\}.
\end{equation}
If the right infimum is attained, the joint minimizers are exactly $(c^\star,q_\theta^\star(c^\star))$ with $c^\star$ minimizing the right side. Compactness of $\cC$ and lower semicontinuity of the reduced objective are sufficient for outer attainment. If $e$, $K_\theta$, and $\cR$ are continuously differentiable near an interior coefficient minimizer, then, for a coefficient direction $h$,
\begin{equation}\label{eq:profile-derivative}
 D_c\Phi_\theta(c)[h]=(D_ce(c)[h])^TA^{-1}e(c)
             -\tfrac12 e(c)^TA^{-1}(D_cK_\theta(c)[h])A^{-1}e(c),
 \quad A=\Sigma+K_\theta(c),
\end{equation}
and $D_c\Phi_\theta+\varrho D\cR=0$ there.
\end{corollary}
\begin{proof}
The inner minimum is attained uniquely for each $c$, so taking coefficient infima preserves equality. The minimizer correspondence follows whenever the outer value is attained. Differentiation uses $D(A^{-1})[h]=-A^{-1}(DA[h])A^{-1}$. Differentiability of $\cR$ is needed only for the stated smooth stationarity condition.
\end{proof}

\subsection{Semidefinite geometry and the factor domain}\label{sec:elimination-factor}
The factor variable below is a deterministic \emph{whitened factor coordinate}, not an identity-covariance Hilbert-valued random element.
\begin{proposition}[Factor form]\label{prop:factor}
Let $C:\cQ\to\cQ$ be bounded, self-adjoint, and positive semidefinite, $T:\cQ\to\R^m$ bounded, and $\Sigma\succ0$. On $\cQ_C:=\Ran(C^{1/2})$, define
\[
 p_C(q):=\min\{\norm{z}_{\cQ}^2:C^{1/2}z=q\}.
\]
Then
\begin{equation}\label{eq:factor}
 \min_{q\in\cQ_C}\{\tfrac12\norm{Tq-e}_{\Sigma^{-1}}^2+\tfrac12p_C(q)\}
 =\min_{z\in\cQ}\{\tfrac12\norm{TC^{1/2}z-e}_{\Sigma^{-1}}^2+\tfrac12\norm{z}^2\}
 =\tfrac12e^T(\Sigma+TCT^*)^{-1}e.
\end{equation}
Both minimizers are unique, and
\[
 z_\star=C^{1/2}T^*(\Sigma+TCT^*)^{-1}e,\qquad
 q_\star=C^{1/2}z_\star=CT^*(\Sigma+TCT^*)^{-1}e.
\]
Moreover,
\begin{equation}\label{eq:factor-stability}
 \norm{z_\star(e)-z_\star(\widetilde e)}
 \le\tfrac12\norm{\Sigma^{-1/2}(e-\widetilde e)}_2.
\end{equation}
\end{proposition}
\begin{proof}
Every preimage of $q\in\cQ_C$ is the minimum-norm preimage $z_0\in(\Ker C)^\perp$ plus an element of $\Ker C$. These kernel components leave $q$ unchanged and strictly increase $\norm{z}^2$ unless zero. Thus minimization over all factors gives the stated source-domain problem, even if the range is not closed in the ambient norm. The factor Hessian is $I+(TC^{1/2})^*\Sigma^{-1}TC^{1/2}\succeq I$, giving a unique $z_\star$ and the formula by the coercive identity with identity geometry. It is orthogonal to $\Ker C$, so it is the minimizing preimage of $q_\star$. Uniqueness of $q_\star$ follows from the strict convexity of the minimum-preimage quadratic. For the last bound set $B=\Sigma^{-1/2}TC^{1/2}$. The factor solution operator is $(I+B^*B)^{-1}B^*\Sigma^{-1/2}$, and its first factor has norm at most $\sup_{s\ge0}s/(1+s^2)=1/2$.
\end{proof}
\begin{remark}[Singular geometry]\label{rem:factor-singular}
Nonuniqueness of a \emph{representation} $q=C^{1/2}z$ is not nonuniqueness of the minimizing factor. If $C\ne0$, then $p_C(q)\ge\norm{q}^2/\norm{C}$; the inverse quadratic retains a strong-convexity lower bound on its admissible linear domain. The issue is that this domain may be proper and the precision unbounded, not disappearance of all curvature. At $q_\star$ both nonnegative energy terms are finite. When $C=0$, the source domain is $\{0\}$ and $z_\star=0$.
\end{remark}

\subsection{Channels and different inferential questions}\label{sec:elimination-readings}
Stacked channels sharing a source are covered by replacing $T$ with $(T_1;\ldots;T_r)$ and using the full measurement weight. The compression has blocks $T_\ell CT_{\ell'}^*$. Shared-source cross-blocks do not in general vanish. Repeated ordinary data vectors, however, do not automatically provide observations of this operator matrix.

\paragraph{Deterministic profile.}
Equation~\eqref{eq:profile} is the value of a penalized source problem. No probability model or determinant term is involved.

\paragraph{Gaussian marginal likelihood.}
Suppose instead that
\[
 d=\mu(c)+T_cq+\varepsilon,\quad q\sim\mathcal N(0,C_{\mathrm{prob}}),\quad
 \varepsilon\sim\mathcal N(0,\Sigma),
\]
independently, with positive trace-class $C_{\mathrm{prob}}$ in an infinite-dimensional source space. Then, for $K_{\mathrm{prob}}(c)=T_cC_{\mathrm{prob}}T_c^*$,
\begin{equation}\label{eq:marginal}
 -\log p(d\mid c)=\tfrac12 e(c)^T(\Sigma+K_{\mathrm{prob}}(c))^{-1}e(c)
       +\tfrac12\log\det(\Sigma+K_{\mathrm{prob}}(c))
       +\tfrac m2\log(2\pi).
\end{equation}
The first term has the profile's algebraic form, with the declared probabilistic compression; equality of the geometries is an additional modeling choice. A constant log determinant is sufficient for the two coefficient objectives to differ only by a constant. It is not necessary for accidental agreement of their minimizers. A nonconstant matrix $K_{\mathrm{prob}}(c)$ need not have nonconstant determinant; orthogonal rotations with isotropic $\Sigma$ give a simple counterexample.

\paragraph{Constrained source fitting.}
For a nonempty admissible set $\cK\subset\cQ$, replacing the penalty by a hard constraint gives
\begin{equation}\label{eq:constrained}
 \inf_{q\in\cK}\tfrac12\norm{T_cq-e(c)}_{\Sigma^{-1}}^2
 =\tfrac12\operatorname{dist}_{\Sigma^{-1}}(e(c),T_c\cK)^2.
\end{equation}
This is a squared-distance value, not generally a support function. If the image is closed and convex, the optimal fitted data are its weighted projection; source uniqueness needs further assumptions. In that case the squared distance is continuously differentiable in $e$, although coefficient regularity needs separate analysis. The quadratic may retain strong convexity in $q$ when $T_c$ is bounded below. We do not analyze this estimator further.

Gaussian posterior optimality, LMMSE optimality under second moments, and GLS/BLUE for deterministic parameters address different statistical questions; see \cref{app:statistics}. Constrained fitting is not a probabilistic reading.

\section{Measurement-relevant geometry}\label{sec:geometry}
\subsection{Structured compression}\label{sec:geometry-compression}
The fixed-Hilbert identity-base family is
\begin{equation}\label{eq:structured}
 C_{\mathrm{pen},\theta}(c)=I+T_c^*S_\theta T_c,\qquad S_\theta\in\PSD^m.
\end{equation}
Writing $G(c)=T_cT_c^*$, its compression is
\begin{equation}\label{eq:structuredK}
 K_\theta(c)=G(c)+G(c)S_\theta G(c).
\end{equation}
It has lower bound $I$ and upper bound $(1+\norm{S_\theta}\norm{T_c}^2)I$, uniformly on bounded families of $S_\theta$ and $T_c$. This is a simple structured family, not a universal parameterization of PSD compressions. In particular $K_\theta\succeq G$. If $G\succ0$, it reaches exactly $\{K:K\succeq G\}$, since $S=G^{-1}(K-G)G^{-1}$. It cannot turn relaxation off by choosing $K=0$ when $G\ne0$.

The reciprocal example will use a positive base $C_0(c)$ rather than $I$. The general base-family formula is
\begin{equation}\label{eq:base-structured}
 C_{\mathrm{pen},\theta}(c)=C_0(c)+C_0(c)T_c^*S_\theta T_cC_0(c),\qquad
 G_0(c)=T_cC_0(c)T_c^*,\quad K_\theta=G_0+G_0S_\theta G_0.
\end{equation}
All quadratic elimination formulas remain pointwise valid under the corresponding coercivity bounds. Calibration theory needs only the specified matrix target $G$; in this base family that target is $G_0$, not $T_cT_c^*$.

\subsection{What compression does and does not identify}\label{sec:geometry-identifiable}
\begin{remark}[Source-geometry non-identifiability]\label{rem:identifiable}
For bounded self-adjoint $C,C'$, equality $T_cCT_c^*=T_cC'T_c^*$ is equivalent to
\[
 \ip{(C-C')T_c^*x}{T_c^*y}_{\cQ}=0\quad(x,y\in\R^m).
\]
Thus only the quadratic form on $\Ran(T_c^*)$ is visible. If $\cQ$ is finite-dimensional and $T_c$ has full column rank, $\Ran(T_c^*)=\cQ$ and this identifies $C$. In an infinite-dimensional source space with finitely many measurements, this cannot occur. For example, any positive operator $N$ supported in $\Ker T_c$ gives the same compression for $C$ and $C+N$, and preserves positivity. Signed perturbations require an explicit positivity check.
\end{remark}
The stated equivalence follows by taking Euclidean inner products of the compressed matrices. Calibration of a compression must therefore not be described as recovery of the physical source covariance.

A separate coefficient-level issue appears in a local linearization
\begin{equation}\label{eq:local}
 \mu(c)\approx\mu(c_0)+J(c_0)(c-c_0),\qquad J(c_0)\in\R^{m\times p}.
\end{equation}
Directions in $\Ker J(c_0)$ remain uninformed by the data for every positive inverse weight. Favorable local risk properties of a matched weight do not imply global nonlinear optimization advantages.

\section{Stability, derivatives, and a reciprocal model}\label{sec:stability}
\subsection{Relative matrix perturbations}\label{sec:stability-perturbation}
\begin{proposition}[Relative perturbation of the quadratic]\label{prop:perturbation}
Let $A\succ0$, let $\Delta$ be symmetric, and suppose $\gamma=\norm{A^{-1/2}\Delta A^{-1/2}}_2<1$. For a residual stack $R\in\R^{m\times r}$, set $\Phi(A;R)=\tfrac12\tr(R^TA^{-1}R)$. Then $A+\Delta\succ0$ and
\begin{equation}\label{eq:perturbation}
 \abs{\Phi(A+\Delta;R)-\Phi(A;R)}
 \le\Phi(A;R)\frac{\gamma}{1-\gamma}.
\end{equation}
\end{proposition}
\begin{proof}
Put $B=A^{-1/2}\Delta A^{-1/2}$ and $Y=A^{-1/2}R$. Then
\[
 \Phi(A+\Delta;R)-\Phi(A;R)
 =-\tfrac12\tr\{Y^T(I+B)^{-1}BY\}.
\]
Since $\norm{(I+B)^{-1}B}_2\le\gamma/(1-\gamma)$, the trace bound gives the result.
\end{proof}
This is an objective-level estimate. Controlling the location of a coefficient minimizer requires additional assumptions; \cref{prop:frozen-stability} supplies one such conditional statement for the frozen estimator.

\subsection{Coefficient regularity and the full product rule}\label{sec:stability-continuity}
\begin{proposition}[Continuity and differentiation]\label{prop:continuity}
Suppose $T_c$, $e(c)$, and $C_{\mathrm{pen},\theta}(c)$ are norm-continuous, with positive geometry and fixed $\Sigma\succ0$. Then $K_\theta(c)$, $(\Sigma+K_\theta(c))^{-1}$, and $\Phi_\theta(c)$ are continuous. If these maps are $C^1$, they are locally Lipschitz and \cref{eq:profile-derivative} holds. In the forcing specialization,
\[
 D_cT_c[h]=-M\cA(c)^{-1}(D_c\cA(c)[h])\cA(c)^{-1}\cB.
\]
For any differentiable geometry in fixed Hilbert structures,
\begin{equation}\label{eq:full-DK}
 DK[h]=(DT[h])CT^*+T(DC[h])T^*+TC(DT[h])^*.
\end{equation}
For \cref{eq:structured}, with $S=S_\theta$ fixed during coefficient differentiation,
\begin{align}
 DC[h]&=(DT[h])^*ST+T^*S\,DT[h],\notag\\
 DG[h]&=DT[h]T^*+T(DT[h])^*,\notag\\
 DK[h]&=DG[h]+DG[h]SG+GS\,DG[h].\label{eq:structured-DK}
\end{align}
If $\Sigma\succeq\sigma_0I$, the inverse satisfies
\[
 \norm{A(c)^{-1}-A(c')^{-1}}_2\le\sigma_0^{-2}\norm{K(c)-K(c')}_2,
 \qquad A(c)=\Sigma+K(c).
\]
It therefore inherits a local Lipschitz bound when $K$ has one, not from continuity alone.
\end{proposition}
\begin{proof}
Products and inversion on uniformly positive matrices preserve continuity. Differentiating $\cA\cA^{-1}=I$ gives the forcing derivative; differentiation of $TCT^*$ and of $G+GSG$ gives the full product rules. The last inequality is the resolvent identity and $\norm{A(c)^{-1}}_2\le\sigma_0^{-1}$.
\end{proof}
For example, $T(c)=c$ and $S=1$ give $K(c)=c^2+c^4$ and $K'(1)=6$, whereas omitting $DC$ gives the incorrect value $4$.

\subsection{An observation-map energy Gram in one dimension}\label{sec:stability-reciprocal}
Let $\Omega=(a,b)\subset\R$, $H=H_0^1(\Omega)$, and choose a \emph{fixed} reference coefficient $c_{\mathrm{ref}}\in L^\infty(\Omega)$ bounded above and below by positive constants. Give $H$ the reference inner product
\[
 \ip{u}{v}_{\mathrm{ref}}=\int_\Omega c_{\mathrm{ref}}u'v'\,dx,
\]
with Riesz map $R_{\mathrm{ref}}:H\to H^*$. For positive $c\in L^\infty(\Omega)$ with uniform ellipticity bounds define
\begin{equation}\label{eq:elliptic1d}
 \ip{\cA(c)u}{v}_{H^*,H}=\int_\Omega c\,u'v'\,dx,\qquad
 P_c:=R_{\mathrm{ref}}^{-1}\cA(c),\qquad C_0(c):=P_c^{-1}.
\end{equation}
Here $P_c$ and $C_0(c)$ are bounded, positive, self-adjoint and uniformly coercive operators on the \emph{fixed} Hilbert space $H$. No equality between an operator inverse and an adjoint of differently typed maps is being asserted.

Let $x_1,\ldots,x_m$ be distinct interior receivers, $Mu=(u(x_i))_{i=1}^m$, and let $M':\R^m\to H^*$ be the dual map, $M'e_i=\delta_{x_i}$. Point evaluation is bounded on $H$ in one dimension. The fixed-reference adjoint is $M^*=R_{\mathrm{ref}}^{-1}M'$. Define $w_i=\cA(c)^{-1}\delta_{x_i}$ and $D_{ij}(c)=w_j(x_i)$.

\begin{lemma}[Reciprocal observation Gram]\label{lem:reciprocity}
In this model,
\begin{equation}\label{eq:DT}
 D(c)=M\cA(c)^{-1}M'=MC_0(c)M^*.
\end{equation}
If $M^{*_{E,c}}$ denotes the adjoint for the energy inner product $\ip{u}{v}_{E,c}=\int c\,u'v'$, then
\begin{equation}\label{eq:gram-reciprocity}
 M^{*_{E,c}}e_i=w_i,\qquad MM^{*_{E,c}}=D(c)=D(c)^T\succeq0.
\end{equation}
For a \emph{state perturbation} $v\in H$, use $T=M$ and the fixed-reference geometry
\begin{equation}\label{eq:reciprocal-geometry}
 C_{\mathrm{pen},\theta}(c)=C_0(c)+C_0(c)M^*S_\theta MC_0(c).
\end{equation}
Its compression is exactly $D(c)+D(c)S_\theta D(c)$.
\end{lemma}
\begin{proof}
Since $C_0(c)M^*=\cA(c)^{-1}R_{\mathrm{ref}}R_{\mathrm{ref}}^{-1}M'=\cA(c)^{-1}M'$, \cref{eq:DT} follows with matching domains. The weak Green equation gives
\[
 \ip{w_i}{v}_{E,c}=v(x_i)\quad(v\in H),\qquad
 \ip{w_i}{w_j}_{E,c}=w_j(x_i)=D_{ij}(c).
\]
This proves the energy-adjoint identity, symmetry, and positivity. Compressing \cref{eq:reciprocal-geometry} with $M$ and $M^*$ gives $D+DS_\theta D$. Coercivity follows from $C_0\succeq\kappa I$ and positivity of the added term; uniform upper bounds follow on bounded $S_\theta$ families.
\end{proof}
The base penalty is $\ip{C_0(c)^{-1}v}{v}_{\mathrm{ref}}=\int c\abs{v'}^2$. At a fixed $c$, \cref{eq:reciprocal-geometry} is precisely the fixed-reference representation of the energy-space family $I+M^{*_{E,c}}S_\theta M$. The fixed-reference representation avoids differentiating a changing Hilbert adjoint implicitly.

\begin{remark}[The forcing map is a different model]\label{rem:reciprocity-scope}
The matrix $D$ is not generally the Gram $TT^*$ of a forcing-to-data operator $T=M\cA(c)^{-1}\cB$. For example, take $(0,1)$, $c=c_{\mathrm{ref}}=1$, one receiver at $1/2$, and the source space $\cQ=H_0^1(0,1)$ with inner product $\ip{q}{v}_{\mathrm{ref}}=\int_0^1q'v'\,dx$. Define $\cB:\cQ\to H^*$ by $\ip{\cB q}{v}_{H^*,H}=\int_0^1qv\,dx$. If $w=\cA(1)^{-1}\delta_{1/2}$, then $Tq=\int_0^1 wq\,dx$ by the weak equation. The adjoint image $z=T^*1$ solves $-z''=w$ with zero boundary values. By symmetry, $w(x)=x/2$ and $z(x)=x/16-x^3/12$ on $[0,1/2]$, with their reflections on $[1/2,1]$. Hence
\[
 TT^*=2\int_0^{1/2} w(x)z(x)\,dx=\tfrac1{480},\qquad D=w(1/2)=\tfrac14.
\]
The adjoint and the numerical forcing Gram depend on the specified source inner product, not on the injection alone. Changing from a forcing perturbation to a state perturbation also changes its penalty geometry.
\end{remark}

\paragraph{Sensitivities.}
In fixed-reference coordinates,
\begin{align}
 DC_0(c)[h]&=-C_0(c)(R_{\mathrm{ref}}^{-1}D\cA(c)[h])C_0(c),\notag\\
 DD(c)[h]&=-M\cA(c)^{-1}(D\cA(c)[h])\cA(c)^{-1}M',\label{eq:DD}\\
 (DD(c)[h])_{ij}&=-\int_\Omega h\,w_i'w_j'\,dx.\notag
\end{align}
The structured derivative is \cref{eq:structured-DK} with $G$ replaced by $D$. These formulas include the derivative of the energy geometry through $C_0(c)$; no such derivative is being discarded.

\paragraph{A matrix-valued acquisition protocol.}
Fix a physical coefficient $c_0$. In replicate session $i$, apply the $m$ unit point forcings $\delta_{x_j}$ one at a time, and record all $m$ receiver responses. If
\[
 Y_i(:,j)=M\cA(c_0)^{-1}\delta_{x_j}+E_i(:,j),\qquad \E E_i=0,
\]
then $X_i=(Y_i+Y_i^T)/2$ is symmetric with mean $D(c_0)$. IID full-session errors with finite second moments give the calibration model below; within-session entrywise independence is unnecessary. This is an idealized one-dimensional protocol with controlled, known excitations. Finite-width sources, excitation errors, or drift require their own mean model. Each replicate uses $m$ excitations and $m^2$ recorded responses. These are \emph{additional measurements}, not free information obtainable by repeating one ordinary data vector. The protocol determines $D(c_0)$, not $D(c)$ for every candidate $c$.

\subsection{Plug-in bias in the structured family}\label{sec:stability-structured}
\begin{proposition}[Structured compression and plug-in bias]\label{prop:matrix}
For $G,S\succeq0$, $K=G+GSG$ has range equal to $\Ran G$. The map $S\mapsto K$ is injective on $\PSD^m$ exactly when $G\succ0$. If a symmetric estimate $\widehat G$ has mean $G$ and finite second moment and $S\succeq0$ is fixed, then
\begin{equation}\label{eq:gram-bias}
 \E[\widehat G+\widehat GS\widehat G]-(G+GSG)
 =\E[(\widehat G-G)S(\widehat G-G)]\succeq0.
\end{equation}
The bias is nonzero exactly when $\E\norm{S^{1/2}(\widehat G-G)}_F^2>0$. Nonzero $S$ and nondegenerate noise alone are not sufficient.
\end{proposition}
\begin{proof}
$K=G^{1/2}(I+G^{1/2}SG^{1/2})G^{1/2}$ has kernel $\Ker G$, hence the same range. Invertible $G$ permits inversion of $GSG$; singular $G$ cannot distinguish $S$ from $S+uu^T$ for $u\in\Ker G$. Expanding $\widehat G=G+(\widehat G-G)$ gives \cref{eq:gram-bias}. Its trace is the stated nonnegative expectation, so it vanishes exactly under that condition.
\end{proof}
If $\widehat G$ is an average of $n$ IID replicates and $S$ is fixed, the bias is $n^{-1}\E[WSW]$. The distinct-replicate construction removes this term without assuming $G$ and $S$ commute. Independently learned random $S$ is handled conditionally on $S$.

\section{Learning and calibration: two distinct reconstruction models}\label{sec:procedure}
\subsection{Exact candidate-dependent relaxation}\label{sec:procedure-exact}
For $\theta=(S_\theta,\varrho_\theta,\Sigma_\theta)$ in a compact admissible parameter set, define
\begin{equation}\label{eq:F-exact}
 F_s^{\mathrm{ex}}(c,\theta)=\tfrac12e_s(c)^T(\Sigma_\theta+K_\theta(c))^{-1}e_s(c)
                             +\varrho_\theta\cR(c).
\end{equation}
In the sample-indexed objective, $e_s(c)=d_s-\mu_s(c)$ uses the fixed training data vector $d_s$ and the modeled mean for its acquisition setup. Here $K_\theta(c)$ is computed from the forward model and the \emph{candidate} coefficient, using either \cref{eq:structuredK} or the declared base family. This is an exact source profile, and \cref{eq:qstar} gives its minimizing latent source. A calibration at one acquisition coefficient is not substituted for $K_\theta(c)$ in this interpretation.

\subsection{Reference-frozen calibrated reconstruction}\label{sec:procedure-frozen}
Let $G_0$ denote the specified Gram at one fixed physical configuration $c_0$, and let $S_\theta$ be learned independently of test calibration replicates. Define $\widehat K_+(S_\theta)$ by \cref{eq:KU,eq:projection}. The frozen estimator is
\begin{equation}\label{eq:F-frozen}
 \widehat c\in\argmin_{c\in\cC_h}\widetilde F(c;\theta),\qquad
 \widetilde F(c;\theta)=\tfrac12e(c)^T(\Sigma_\theta+\widehat K_+(S_\theta))^{-1}e(c)
                         +\varrho_\theta\cR(c).
\end{equation}
The matrix in this objective is fixed during optimization over $c$. The oracle version replaces $\widehat K_+$ by $G_0+G_0S_\theta G_0$. Neither is generally the exact profile \cref{eq:F-exact}. In particular the coefficient gradient has \emph{no} $D_cK$ term for the frozen matrix.

One must not replace only the denominator of \cref{eq:qstar} by $\Sigma+\widehat K_+$ and call the result an exact minimizing source. For example, $T=C=\Sigma=e=1$ has exact minimizer $q=1/2$, whereas this substitution with $\widehat K_+=0$ gives $q=1$, whose source gradient is nonzero. The frozen procedure reports a coefficient and an inverse weight, not an exact relaxed source. If a source estimate is subsequently desired, one must specify and solve a separate source problem. PSD projection alone neither constructs a source-space lift nor preserves the family $G+GSG$.

\subsection{Matched training and independent test calibration}\label{sec:procedure-steps}
An operational two-stage procedure is as follows.
\begin{enumerate}[label=(\roman*),leftmargin=*]
\item Specify the forcing or state-perturbation model, the fixed Hilbert geometry, $\cC_h\subset\R^p$, and $\Theta$. Enforce $S_\theta\succeq0$, $\varrho_\theta\ge0$, and $\Sigma_\theta\succeq\sigma_0 I$ on $\Theta$.
\item For each sample $s$, assume that its reference coefficient $c_s^{\mathrm{true}}$ is known for supervision (for example, in simulation or a controlled training experiment) and that the data vector $d_s$ is acquired at this configuration. The pair $(c_s^{\mathrm{true}},d_s)$ is fixed during training. The lower-level candidate $c$ and reconstruction $\widehat c_s$ are distinct from this reference. Choose the deployed estimator: exact \cref{eq:F-exact}, or frozen \cref{eq:F-frozen}. For matched frozen training, also collect $n_s\ge2$ matrix replicates $X_{s,i}$ at the same physical configuration $c_s^{\mathrm{true}}$, with mean $G(c_s^{\mathrm{true}})$ in the declared calibration model. The replicates are fixed inputs to training. Set
\begin{equation}\label{eq:F-frozen-training}
 F_s^{\mathrm{fr}}(c,\theta)=\tfrac12 e_s(c)^T
       (\Sigma_\theta+\widehat K_{s,+}(S_\theta))^{-1}e_s(c)
       +\varrho_\theta\cR(c).
\end{equation}
Here $\widehat K_{s,+}(S_\theta)$ is formed from the stored training replicates but is \emph{independent of the candidate} $c$. Oracle frozen training using a computed $G(c_s^{\mathrm{true}})$ is an alternative idealization, not the same finite-replicate estimator. In both acquisition-based lower objectives, the reference coefficient is not passed in as a reconstruction input; it labels the acquisition configuration and supplies the upper-level target.
\item Minimize the supervised upper loss in \cref{eq:bilevel} for the chosen $F_s$. Smooth-branch hypergradients apply only when the relevant assumptions hold. For training with spectral projection, a sufficient local smoothness condition is that the raw estimate's spectrum stays away from zero, together with smooth parameterization of $S_\theta$ and $\Sigma_\theta$. At zero-eigenvalue crossings, the smooth formulas are not justified.
\item Keep the learned $\widehat\theta$ fixed. For frozen deployment, acquire fresh IID matrix replicates at the test physical configuration, independent of \emph{all} training data and the fitted parameter. Form $\widehat K_U$ and apply the PSD admissibility rule. Reconstruct with \cref{eq:F-frozen}. For exact deployment, continue to evaluate the candidate-dependent \cref{eq:F-exact} instead.
\end{enumerate}
Each training or test matrix-replicate session has the extra measurement cost stated in the reciprocal protocol. Comparisons against methods without these acquisitions must disclose this data-budget difference.

\subsection{The role of independence}\label{sec:procedure-independent}
The raw calibration theorem applies at a fixed $S$, or conditionally on an independently learned $S$. It does not assert unbiasedness after optimizing $S$ on the \emph{same} calibration replicates. Dependence of $S$ on training replicates is normal in supervised learning; those training data are not used for the independent test calibration guarantee. If the test data vector and test calibration are also independent, one may condition on the data vector in the reconstruction bound below. Reusing test observations requires a separate conditional analysis.

Training with \cref{eq:F-exact} and testing with \cref{eq:F-frozen} is a change of estimator. It can be explored empirically, but the exact-profile identities do not justify its performance. Similarly, estimating $G(c_0)$ provides no estimate of the whole function $c\mapsto G(c)$ without further modeling and acquisitions.

\section{Bilevel existence and local hypergradients}\label{sec:bilevel}
\subsection{Compactness before sensitivity}\label{sec:bilevel-existence}
For either matched lower model above, let $\ell_s(c;c_s^{\mathrm{true}})$ be a real-valued loss comparing a reconstruction $c\in\cC_h$ to the known training reference. For instance, when the two coefficients are represented on the same grid, use $\ell_s(c;c_s^{\mathrm{true}})=\tfrac12\norm{c-c_s^{\mathrm{true}}}_2^2$. The references are held fixed, and the supervised training problem is
\begin{equation}\label{eq:bilevel}
 \begin{aligned}
 \min_{\substack{\theta\in\Theta,\;(\widehat c_1,\ldots,\widehat c_{n_{\mathrm{tr}}})
                  \in\cC_h^{n_{\mathrm{tr}}}}}\quad
 &\frac1{n_{\mathrm{tr}}}\sum_{s=1}^{n_{\mathrm{tr}}}
      \ell_s(\widehat c_s;c_s^{\mathrm{true}})+\Psi_0(\theta)\\
 \text{subject to}\quad
 &\widehat c_s\in\argmin_{c\in\cC_h}F_s(c,\theta),\qquad s=1,\ldots,n_{\mathrm{tr}}.
 \end{aligned}
\end{equation}
When the lower minimizer is unique, write it as $\widehat c_s(\theta)$ and substitute it in the upper loss. Otherwise \cref{eq:bilevel} is optimistic: the upper problem can select among global lower minimizers. Local numerical solves need not find these global minimizers. The supervised formulation assumes known references for training, not for test reconstruction.

\begin{assumption}[Finite-dimensional compact bilevel setting]\label{ass:bilevel}
The sets $\cC_h\subset\R^p$ and $\Theta\subset\R^{d_\theta}$ are nonempty and compact, each $F_s$ is finite and jointly continuous on $\cC_h\times\Theta$, and each loss $\ell_s(\cdot;c_s^{\mathrm{true}})$ is continuous on $\cC_h$. The upper regularizer $\Psi_0:\Theta\to\R\cup\{+\infty\}$ is lower semicontinuous, bounded below, and finite at at least one parameter. In particular any parameter-dependent $\varrho_\theta\cR(c)$ is included in the full continuity requirement on $F_s$.
\end{assumption}

\begin{proposition}[Unique lower solutions when strongly convex]\label{prop:unique-bilevel}
Under \cref{ass:bilevel}, suppose additionally that $\cC_h$ is convex and every $F_s(\cdot,\theta)$ is $\mu_{\mathrm{sc}}$-strongly convex with a common $\mu_{\mathrm{sc}}>0$. Each lower solution $\widehat c_s(\theta)$ is unique and continuous. The reduced upper loss is lower semicontinuous and attains its minimum. It is continuous if $\Psi_0$ is finite and continuous.
\end{proposition}
\begin{proof}
Existence follows from compactness and uniqueness from strong convexity. For $x=\widehat c_s(\theta)$ and $x'=\widehat c_s(\theta')$, optimality and strong convexity give
\[
 \tfrac{\mu_{\mathrm{sc}}}{2}\norm{x'-x}^2
 \le F_s(x',\theta)-F_s(x,\theta)
 \le 2\sup_{c\in\cC_h}\abs{F_s(c,\theta)-F_s(c,\theta')}.
\]
The supremum tends to zero by uniform continuity on the compact product. The upper claims follow by composition and compactness. This comparison uses the full lower objective, including the variation of $\varrho_\theta$.
\end{proof}

\begin{proposition}[Genuinely optimistic attainment]\label{prop:optimistic}
Under \cref{ass:bilevel}, without lower uniqueness, problem \cref{eq:bilevel} has a minimizer.
\end{proposition}
\begin{proof}
Each lower argmin is nonempty. Its graph is closed: if $\theta_k\to\theta$ and $c_k\to c$ with $F_s(c_k,\theta_k)\le F_s(v,\theta_k)$ for every $v\in\cC_h$, continuity gives the same inequality at $(c,\theta)$. Consequently the joint feasible graph for all training samples is a nonempty closed subset of the compact product $\Theta\times\cC_h^{n_{\mathrm{tr}}}$. The upper objective is lower semicontinuous on that graph and finite somewhere, so it attains its minimum.
\end{proof}
These are existence results, not Lipschitz or Clarke-stationarity results. Even uniform strong convexity and joint continuity can give a non-Lipschitz solution map: $F(c,\theta)=(c-\sqrt{|\theta|})^2$ on $[0,1]\times[-1,1]$ has solution $\sqrt{|\theta|}$. No Clarke subgradient of the reduced loss is invoked under these assumptions. Strong convexity here is an extra coefficient-level hypothesis, not a consequence of source elimination.

\subsection{Interior branches}\label{sec:bilevel-interior}
\begin{proposition}[Interior hypergradient]\label{prop:interior}
Let $F(c,\theta)$ be $C^2$ near an interior pair $(c_0,\theta_0)$ with $\nabla_cF(c_0,\theta_0)=0$ and $H_0=\nabla^2_{cc}F(c_0,\theta_0)\succ0$. There is a unique local $C^1$ stationary branch $c^\star(\theta)$ near that pair, consisting of strict local minimizers, with
\begin{equation}\label{eq:hypergrad-interior}
 D_\theta c^\star=-H^{-1}F_{c\theta},\qquad
 H=F_{cc}(c^\star(\theta),\theta).
\end{equation}
For a $C^1$ loss $\ell(c,\theta)$,
\begin{equation}\label{eq:hypergrad-loss}
 D_\theta\ell(c^\star(\theta),\theta)
 =\ell_\theta-\nabla_c\ell^TH^{-1}F_{c\theta}.
\end{equation}
This is a local branch formula; it differentiates the global reduced upper objective only if that branch is the selected global lower solution.
\end{proposition}
\begin{proof}
Apply the implicit-function theorem to $\nabla_cF=0$. Positive definiteness persists locally. Differentiating stationarity and then the loss gives the two formulas.
\end{proof}
An adjoint solve $H^Ta=\nabla_c\ell$ computes the whole upper gradient without one solve per parameter component. This says nothing about the conditioning or cost of the PDE and Hessian solves. For the supervised loss, use $\ell(c,\theta)=\ell_s(c;c_s^{\mathrm{true}})$ in the local branch formulas: its explicit $\theta$-derivative is zero because the reference and loss are fixed. Add $D\Psi_0$ and average the sample contributions when assembling a smooth upper gradient.

\subsection{Smooth fixed-constraint active sets}\label{sec:bilevel-active}
The next statement concerns smooth objectives with smooth constraints. It does not directly cover TV, TGV, or $\ell^1$ penalties. Such penalties require a justified smooth reformulation/active-manifold analysis or a different nonsmooth sensitivity theory.

\begin{proposition}[Active-set hypergradient]\label{prop:active}
Let $F$ be $C^2$ near $(c_0,\theta_0)$ and let the \emph{parameter-independent} constraints be $g_a(c)=0$ and $g_b(c)\le0$, with $C^2$ constraint functions. Suppose $c_0$ is a local minimizer satisfying LICQ, strict complementarity on its active inequality set $I$, and the strong second-order sufficient condition. With the active constraints stacked as $g(c)$, multipliers $\nu$, and Lagrangian $L(c,\theta,\nu)=F(c,\theta)+\nu^Tg(c)$, there is a local $C^1$ primal--dual branch with the same active set. Define
\[
 \mathcal K=\begin{pmatrix}L_{cc}&g_c^T\\g_c&0\end{pmatrix}.
\]
Then
\begin{equation}\label{eq:KKT-derivative}
 \mathcal K\begin{pmatrix}D_\theta c^\star\\D_\theta\nu^\star\end{pmatrix}
 =-\begin{pmatrix}L_{c\theta}\\0\end{pmatrix}.
\end{equation}
For a $C^1$ loss $\ell(c,\theta)$,
\begin{equation}\label{eq:hypergrad-active}
 D_\theta\ell(c^\star(\theta),\theta)=\ell_\theta-\alpha_\ell^TL_{c\theta},\qquad
 \mathcal K^T\begin{pmatrix}\alpha_\ell\\\beta_\ell\end{pmatrix}
 =\begin{pmatrix}\nabla_c\ell\\0\end{pmatrix}.
\end{equation}
\end{proposition}
\begin{proof}
LICQ and the positive reduced Lagrangian Hessian make the active KKT matrix nonsingular. The implicit-function theorem gives a $C^1$ branch. Strict complementarity and strict inactivity of the remaining constraints preserve the active set locally; the second-order condition preserves strict local minimality. Differentiating the KKT equations gives \cref{eq:KKT-derivative}, with zero lower right-hand side because $g$ is independent of $\theta$. Multiplication by the adjoint solution and the chain rule give \cref{eq:hypergrad-active}; see also \cite{BonnansShapiro2000,DontchevRockafellar2009}.
\end{proof}
Neither branch formula covers active-set transitions or selection among competing local minima.

\subsection{Derivatives for the two deployed models}\label{sec:bilevel-reduced}
For exact relaxation, use \cref{eq:profile-derivative,eq:full-DK,eq:structured-DK}. At fixed $c$,
\[
 D_SK(c)[\delta S]=G(c)\delta S G(c),\qquad
 D_cD_SK(c)[h,\delta S]=DG(c)[h]\delta S G(c)+G(c)\delta S DG(c)[h].
\]
The same formulas hold with the declared base Gram $G_0$ or $D$. In the reciprocal model, $DD$ is \cref{eq:DD}; ignoring the derivative of $C_0$ would be incorrect.

For frozen reconstruction, $D_c\widehat K_+=0$, so
\[
 D_c\widetilde\Phi(c)[h]=(D_ce(c)[h])^T
                    (\Sigma_\theta+\widehat K_+(S_\theta))^{-1}e(c).
\]
Upper-parameter derivatives still include the dependence of $\Sigma_\theta$, $S_\theta$, and $\varrho_\theta$, and, for matched finite-replicate training, the spectral projection where it is differentiable. Candidate-independent does not mean parameter-independent.

\section{Finite-replicate calibration and its reconstruction interface}\label{sec:calibration}
\subsection{An IID matrix model and exact variance}\label{sec:calibration-estimator}
Fix an acquisition configuration. Let its specified matrix target be $G\succeq0$. In the reciprocal protocol $G=D(c_0)$; in any other model a mean-$G$ matrix acquisition protocol is an additional assumption, not a consequence of repeated scalar/vector data. For $n\ge2$, define
\begin{equation}\label{eq:KU}
 \widehat K_U=\bar X+\frac1{n(n-1)}\sum_{i\ne j}X_iSX_j,\qquad
 \bar X=\frac1n\sum_iX_i.
\end{equation}
The sum uses ordered pairs. Since the $X_i$ and $S$ are symmetric, pairing $(i,j)$ and $(j,i)$ makes the result symmetric even when these matrices do not commute.

For a square-integrable matrix $Z$, define the \emph{scalar} Frobenius variance
\[
 \Var_F(Z\mid S):=\E[\norm{Z-\E[Z\mid S]}_F^2\mid S].
\]
This is the trace of the conditional covariance of its vectorization, not an entrywise matrix variance.

\begin{theorem}[Conditional unbiasedness and exact IID variance]\label{thm:calibration}
Let $X_1,\ldots,X_n$ be IID symmetric random matrices with common mean $G\succeq0$ and $\E\norm{X_1-G}_F^2<\infty$. Let $S\succeq0$ be deterministic, or an almost surely finite symmetric PSD random matrix independent of the full replicate collection. Put $K(S)=G+GSG$. Then
\begin{equation}\label{eq:calibration-unbiased}
 \E[\widehat K_U\mid S]=K(S),
\end{equation}
and, for $W=X_1-G$ and an independent copy $W'$,
\begin{align}
 \Var_F(\widehat K_U\mid S)&=\frac{a(S)}n+\frac{2b(S)}{n(n-1)}=:v_n(S),\label{eq:calibration-variance}\\
 a(S)&=\E[\norm{\mathcal L_S(W)}_F^2\mid S],\qquad
 b(S)=\E[\norm{\mathcal H_S(W,W')} _F^2\mid S],\notag\\
 \mathcal L_S(A)&=A+GSA+ASG,\qquad
 \mathcal H_S(A,B)=\tfrac12(ASB+BSA).\notag
\end{align}
These conditional quantities are finite for almost every fixed $S$. Unconditional bounds for random $S$ additionally require integrability of $a(S)$ and $b(S)$, for example bounded $S$.
\end{theorem}
The proof is in \cref{app:calibration}. This is an exact finite-sample order-two U-statistic variance decomposition. Only second moments of individual replicates are needed because products involve \emph{independent distinct} replicates. The nonnegative pair contribution is a cost of unbiasedly estimating the quadratic target; no variance superiority over $\bar X$ or other estimators is asserted, and $\bar X$ generally estimates a different target.

\subsection{PSD projection is a sufficient admissibility rule}\label{sec:calibration-admissible}
\begin{proposition}[Admissible projection and possible bias]\label{prop:admissible}
For symmetric $A=\sum_k\rho_kv_kv_k^T$, define
\begin{equation}\label{eq:projection}
 \Pi_{\PSD^m}(A)=\sum_k\max(\rho_k,0)v_kv_k^T,\qquad
 \widehat K_+=\Pi_{\PSD^m}(\widehat K_U).
\end{equation}
Then $\Sigma+\widehat K_+\succeq\Sigma\succ0$, and for every $K\succeq0$,
\begin{equation}\label{eq:admissible}
 \norm{\widehat K_+-K}_F\le\norm{\widehat K_U-K}_F.
\end{equation}
The projected estimator need not be unbiased; this failure occurs under the IID model itself.
\end{proposition}
\begin{proof}
The PSD cone is closed and convex in the Frobenius Hilbert space, and the spectral formula is its unique metric projection. Nonexpansiveness, together with $\Pi_{\PSD^m}(K)=K$, gives \cref{eq:admissible}. The following IID example proves possible bias; nonlinearity of a projection alone would not prove it.
\end{proof}
\begin{remark}[An IID projection-bias example]\label{rem:counterexample}
Let $m=1$, $G=S=1$, $n=2$, and let $X_1,X_2$ be IID with probabilities $1/2$ at $3$ and $-1$. They have mean $1$. The raw estimator $\widehat K_U=(X_1+X_2)/2+X_1X_2$ takes values $12,-2,0$ with probabilities $1/4,1/2,1/4$. Therefore
\[
 \E\widehat K_U=2=G+GSG,\qquad \E\widehat K_+=3\ne2.
\]
On the negative event, $\Sigma=1$ gives $\Sigma+\widehat K_U=-1$, which is not a positive quadratic weight. Projection gives $0$ for that realization. This separates realized indefiniteness from bias in expectation.
\end{remark}
PSD projection is not the only way to ensure a positive inverse weight. A slightly negative raw compression can still leave $\Sigma+\widehat K_U\succ0$. Nor does PSD projection enforce $\widehat K_+\succeq G$, the range constraint of a given $T_c$, or membership in $G+GSG$. It is used here as a sufficient measurement-space admissibility device for the frozen estimator.

\subsection{Calibration and inverse-weight consistency}\label{sec:calibration-consistency}
\begin{corollary}[Mean-square calibration and inverse-weight bounds]\label{cor:calibration-consistency}
Under \cref{thm:calibration},
\begin{equation}\label{eq:calibration-L2}
 \E[\norm{\widehat K_+-K(S)}_F^2\mid S]\le v_n(S).
\end{equation}
If $\Sigma\succeq\sigma_0 I$, with $\sigma_0>0$, is fixed or independently learned along with $S$, then
\begin{equation}\label{eq:inverse-L2}
 \E[\norm{(\Sigma+\widehat K_+)^{-1}-(\Sigma+K(S))^{-1}}_F^2\mid S,\Sigma]
 \le\sigma_0^{-4}v_n(S).
\end{equation}
For every fixed $S$, both estimators are mean-square consistent, with root-mean-square error $O(n^{-1/2})$ (or faster in degenerate cases), and hence consistent in probability.
\end{corollary}
\begin{proof}
Square \cref{eq:admissible}, condition on $S$, and use unbiasedness and \cref{eq:calibration-variance}. The resolvent identity gives
\[
 \norm{\widehat A^{-1}-A^{-1}}_F
 =\norm{\widehat A^{-1}(K-\widehat K_+)A^{-1}}_F
 \le\sigma_0^{-2}\norm{\widehat K_+-K}_F,
\]
where $\widehat A=\Sigma+\widehat K_+$ and $A=\Sigma+K$. Independence of the jointly learned parameters from the calibration collection preserves the conditional replicate law. Squaring gives the second bound. Markov's inequality implies convergence in probability.
\end{proof}
The bounds concern a fixed acquisition target. They do not imply population consistency of supervised learning, unbiased finite-sample inverse weights, or equality of frozen and candidate-dependent profiles. Crude tail bounds follow from Markov's inequality, but no sharp confidence calibration or asymptotic distribution is supplied.

\subsection{A conditional link to the frozen coefficient estimator}\label{sec:calibration-reconstruction}
The next statement compares estimated and \emph{oracle frozen} weights. The oracle coefficient need not be the true physical coefficient.
\begin{proposition}[Frozen-weight solution stability]\label{prop:frozen-stability}
Let $\cC_h\subset\R^p$ be nonempty, compact, and convex. Let $e$ and $\cR$ be $C^1$ on a neighborhood, and suppose
\[
 \sup_{c\in\cC_h}\norm{e(c)}_2\le E_0,\qquad
 \sup_{c\in\cC_h}\norm{De(c)}_2\le E_1.
\]
For $K,\widehat K\succeq0$, $\Sigma\succeq\sigma_0 I$, and fixed $\varrho\ge0$, define
\begin{align*}
 F_A(c)&=\tfrac12e(c)^TA^{-1}e(c)+\varrho\cR(c), & A&=\Sigma+K,\\
 F_{\widehat A}(c)&=\tfrac12e(c)^T\widehat A^{-1}e(c)+\varrho\cR(c), & \widehat A&=\Sigma+\widehat K.
\end{align*}
Assume the oracle objective $F_A$ is $\mu_{\mathrm{sc}}$-strongly convex on $\cC_h$. Its minimizer $c_A$ is unique. Every global minimizer $c_{\widehat A}$ of $F_{\widehat A}$ on this set satisfies
\begin{equation}\label{eq:frozen-stability}
 \norm{c_{\widehat A}-c_A}_2
 \le\frac{E_0E_1}{\mu_{\mathrm{sc}}}\norm{\widehat A^{-1}-A^{-1}}_2
 \le\frac{E_0E_1}{\mu_{\mathrm{sc}}\sigma_0^2}\norm{\widehat K-K}_F.
\end{equation}
For the independently calibrated $\widehat K_+$ and $K=K(S)$, condition on the learned parameters and, if random, an independent test data vector. Then
\begin{equation}\label{eq:coefficient-L2}
 \E[\norm{c_{\widehat A}-c_A}_2^2\mid S,\Sigma,e]
 \le\left(\frac{E_0E_1}{\mu_{\mathrm{sc}}\sigma_0^2}\right)^2v_n(S),
\end{equation}
where constants may depend on the conditioned objects. Any random minimizer selection in this statement is assumed measurable.
\end{proposition}
\begin{proof}
Compactness gives attainment. Put $u=c_{\widehat A}-c_A$. The first-order variational inequalities at the two constrained minimizers and strong monotonicity of $\nabla F_A$ give
\[
 \mu_{\mathrm{sc}}\norm{u}^2\le\ip{\nabla F_A(c_{\widehat A})-\nabla F_A(c_A)}{u}
 \le\ip{\nabla F_A(c_{\widehat A})-\nabla F_{\widehat A}(c_{\widehat A})}{u}.
\]
The gradient difference is $De(c)^T(A^{-1}-\widehat A^{-1})e(c)$, whose norm is bounded by $E_0E_1\norm{A^{-1}-\widehat A^{-1}}_2$. Cancel $\norm{u}$ when nonzero and use the resolvent bound. Squaring and applying \cref{cor:calibration-consistency}, with the stated independence, proves \cref{eq:coefficient-L2}.
\end{proof}
Strong convexity of the estimated objective was not required, but a global minimizer on the stated set was. For a local application, take a compact convex neighborhood where oracle strong convexity holds and minimize over that neighborhood; this gives no guarantee that a global nonlinear solver enters it. The bound isolates calibration error from oracle bias and optimization error, neither of which it controls.

\begin{remark}[When does an oracle advantage survive calibration?]\label{rem:oracle-advantage} 
Let $c_A$, $c_{\widehat A}$, and $c_{\mathrm{base}}$ denote the oracle frozen,
calibrated frozen, and baseline reconstructions. Fix the independently
learned parameters, if any, and work under the joint law of the test instance
and the fresh calibration replicates. Conditional on the physical configuration,
assume that the calibration replicates are IID under the corresponding matrix
model and independent of the test reconstruction data. Assume that all three
reconstruction errors are square-integrable. For each test instance, $c_A$
and $c_{\widehat A}$ use the same residual, feasible set, and regularization
term and differ only in their frozen weights.

Define $r(z)=(\E\norm{z-c^{\mathrm{true}}}_2^2)^{1/2}$, with expectation under
the joint test and calibration law. For reconstructions that are functions
of the test variables alone, this reduces to expectation over the marginal
test law. If $B_n<\infty$ satisfies
$(\E\norm{c_{\widehat A}-c_A}_2^2)^{1/2}\le B_n$, the triangle inequality in
$L^2$ gives $r(c_{\widehat A})\le r(c_A)+B_n$. Thus an oracle margin
$r(c_{\mathrm{base}})-r(c_A)>B_n$ is sufficient for the calibrated estimator
to have strictly smaller test RMSE than the baseline.

Under the hypotheses of \cref{prop:frozen-stability}, apply
\cref{eq:coefficient-L2} conditionally on the physical configuration and test
reconstruction data, then average over their marginal law to obtain
\[
B_n^2=\mathbb E_{\mathrm{test}}\!\left[
\left(\frac{E_0 E_1}{\mu_{\mathrm{sc}}\sigma_0^2}\right)^2 v_n(S)\right],
\]
provided the constants are measurable with respect to the conditioned
variables and the right-hand side is finite. When the test law includes
multiple physical configurations, evaluate $v_n(S)$ and the stability
constants for the corresponding configuration and acquisition protocol.

The criterion neither establishes an oracle advantage nor guarantees that a
nonlinear solver finds the minimizers. If the oracle margin does not exceed
$B_n$, the criterion is inconclusive. This does not imply that the calibrated
estimator is worse. Comparisons should also account for the additional
calibration data.
\end{remark}

\paragraph{Elliptic example: from noisy responses to a frozen weight.}
Consider the one-dimensional Dirichlet problem with colocated point
sources and receivers, as in Section~3.1.2 of \cite{vanLeeuwenYang2025}.
On $(0,1)$, take the calibration coefficient $c_0$ equal to $1$ on
$(0,1/2)$ and $2$ on $(1/2,1)$, with $m=8$ locations $x_j=j/(m+1)$.
For a unit point source at $x_j$, let $u_j$ solve
$-(c_0 u_j')'=\delta_{x_j}$ with zero boundary values. Writing
$r(x)=\int_0^x c_0(t)^{-1}\,dt$, the response and energy Gram are
\[
 (G_0)_{ij}=u_j(x_i)=\int_0^1 c_0 u_i' u_j'\,dx
 =r(\min\{x_i,x_j\})-\frac{r(x_i)r(x_j)}{r(1)}.
\]
Indeed, the Green function on the right vanishes at the endpoints, is
continuous, and has a flux jump of $-1$ across the source. This exact
response avoids PDE discretization error, isolating the effect of noisy
calibration in this synthetic setting.

\begin{figure}[tbp]
 \centering
 \includegraphics[width=.94\linewidth]{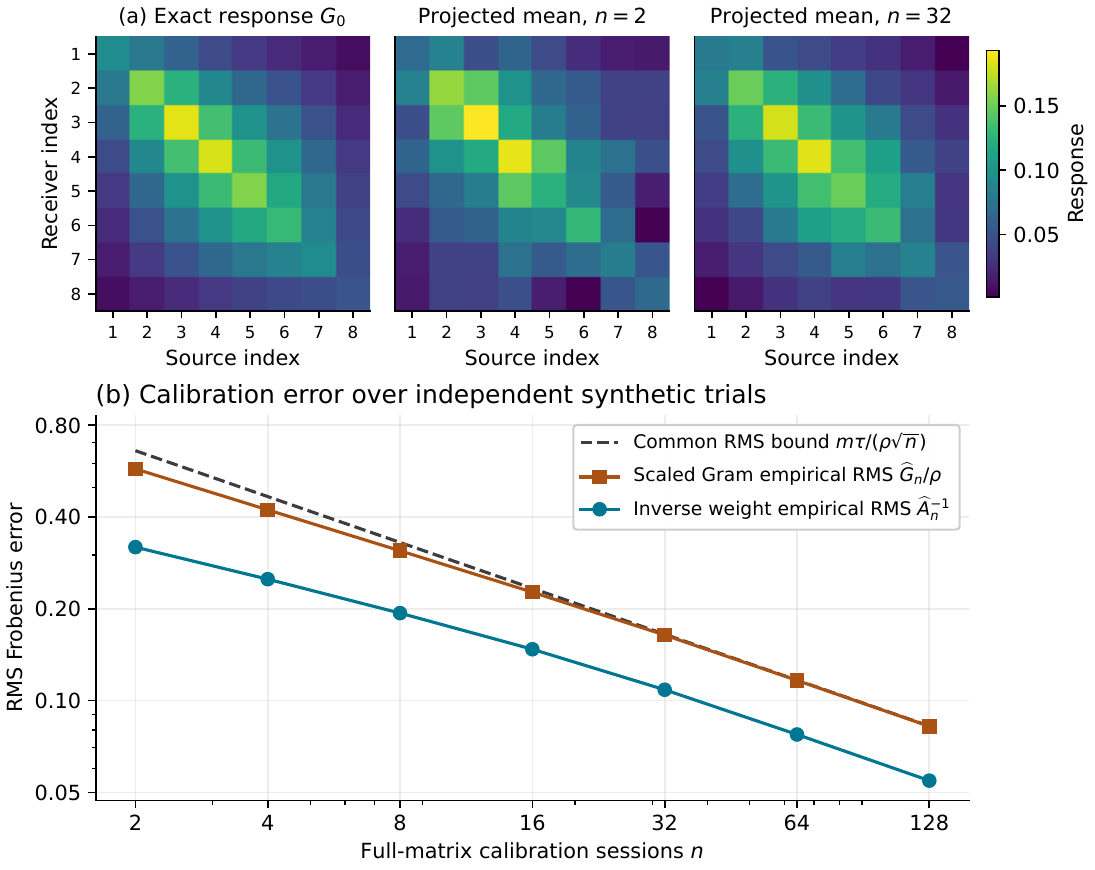}
 \caption{Synthetic one-dimensional elliptic calibration. (a) Exact
 response and PSD-projected averages from the first $2$ and all $32$
 sessions of one illustrative draw on a common color scale. (b) Empirical
 RMS Frobenius errors of the scaled Gram $\widehat G_n/\rho$ and the inverse
 weight $\widehat A_n^{-1}$ relative to their exact counterparts. The dashed
 curve is the common theoretical RMS upper bound for both errors, not a
 per-trial threshold. Each synthetic session represents eight source
 excitations with eight receiver readings per excitation, modeled as a
 symmetric matrix. No reconstruction data are used.}
 \label{fig:elliptic-calibration}
\end{figure}

A synthetic calibration session records a full matrix $X_\ell=G_0+Z_\ell$.
The $Z_\ell$ are IID symmetric matrices, with independent
$\mathcal N(0,\tau^2)$ upper-triangular entries (including the diagonal)
and $\tau=0.035$. Equivalently, the noise for each unordered pair of
locations is drawn once per session and mirrored in the matrix. This
idealized observation model enforces reciprocity; it does not model
independent noise in the two reciprocal readings. For $\rho=0.30$, form
\[
 \widehat G_n=\Pi_{\PSD^m}\!\left(\frac1n\sum_{\ell=1}^n X_\ell\right),
 \qquad
 \widehat A_n^{-1}=(I+\rho^{-1}\widehat G_n)^{-1},
 \qquad A^{-1}=(I+\rho^{-1}G_0)^{-1}.
\]
Only the session matrices enter the estimated weight; the exact Gram is
used for scoring. Applying \cref{eq:inverse-L2} to $\rho^{-1}X_\ell$
with $S=0$, $\Sigma=I$, and $\sigma_0=1$ gives the bound
\[
 \left(\mathbb E\|\widehat A_n^{-1}-A^{-1}\|_F^2\right)^{1/2}
 \le \frac{m\tau}{\rho\sqrt n}.
\]
Here $\mathbb E\|Z_\ell\|_F^2=m^2\tau^2$ (counting each off-diagonal
pair twice), so the scaled sample mean has Frobenius mean-square error
$m^2\tau^2/(\rho^2 n)$. By nonexpansiveness of the PSD projection and
$G_0\succeq0$, the same RMS bound controls
$\bigl(\mathbb E\|\widehat G_n/\rho-G_0/\rho\|_F^2\bigr)^{1/2}$.
\Cref{fig:elliptic-calibration} compares the common bound with the empirical RMS
Frobenius error over $4000$ independent trials per $n$. For example, the
inverse-weight error decreases from $0.319$ at $n=2$ to $0.109$ at $n=32$,
versus bounds of $0.660$ and $0.165$, respectively. The experiment
illustrates calibration of the frozen weight. It does not test
reconstruction error, retention of an oracle RMSE advantage
(\cref{rem:oracle-advantage}), or preservation of objective-landscape
properties. Dependent noise, session drift, and source calibration errors
are outside this experiment.

\paragraph{Seismic-motivated RMSE-margin illustration.}
Section~3.5 of \cite{vanLeeuwenYang2025} compares conventional least squares,
a candidate-dependent metric, and a frozen data-driven metric for
velocity inversion with the two-dimensional Helmholtz equation. The reported
reconstructions motivate a separate question: if an exact frozen metric
improves velocity RMSE, is the improvement large enough to survive noisy
calibration? The published frozen data-driven metric is not the exact
frozen-weight oracle of this comparison.

\begin{figure}[tbp]
 \centering
 \includegraphics[width=.98\linewidth]{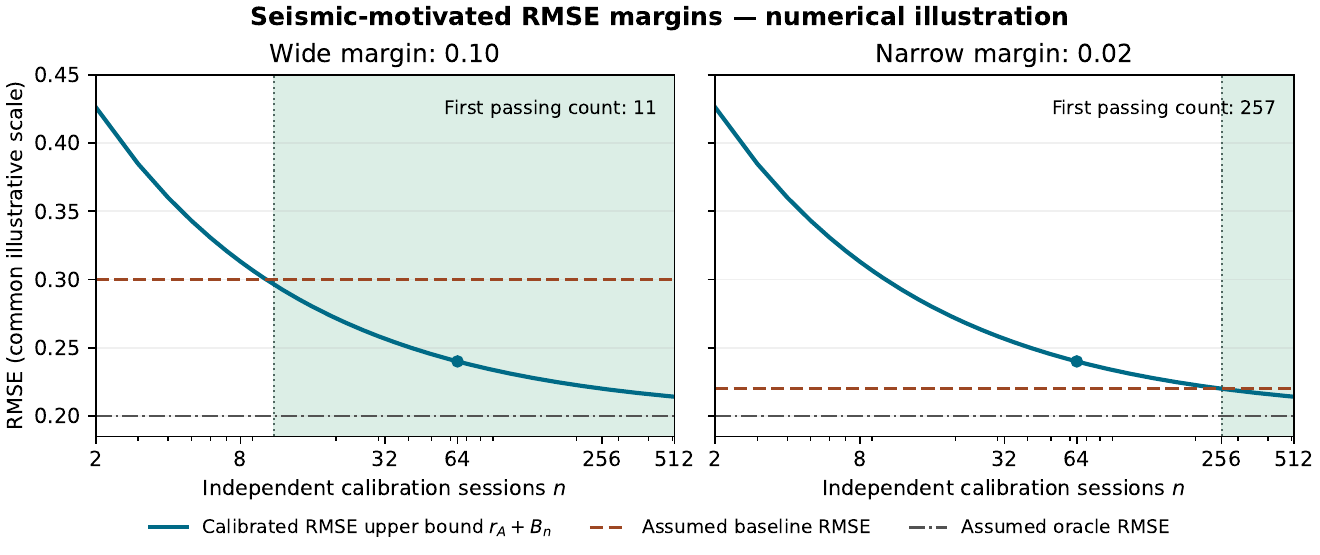}
 \caption{Seismic-motivated numerical illustration, not an inversion
 experiment. Both scenarios use the same assumed oracle RMSE and calibration
 bound. Under these assumptions, the shaded region marks integer session
 counts for which the upper bound is strictly below the baseline. Dots
 mark $n=64$. Outside the shaded region the sufficient condition is
 inconclusive. The curve is an upper bound, not a measured error.}
 \label{fig:seismic-margin}
\end{figure}

In the proposed comparison, take conventional least squares as the baseline
and compare it with an exact frozen-weight oracle and an independently
calibrated frozen-weight estimator. Use the same test data, feasible set,
and regularization in all three reconstructions. In particular, only the
frozen weight changes between the oracle and calibrated reconstructions.
Write $r_A=r(c_A)$ and $r_{\mathrm{base}}=r(c_{\mathrm{base}})$ for the
velocity RMSEs in the common scale of \cref{rem:oracle-advantage}.
The calibrated RMSE and the displacement bound are
understood under the joint test and calibration law,
using the same coefficient norm and normalization.
If $B_n$ is a valid bound on the RMS displacement from the oracle, then
\[
 r(c_{\widehat A})\le r_A+B_n,
 \qquad
 B_n<r_{\mathrm{base}}-r_A
 \quad\Longrightarrow\quad r(c_{\widehat A})<r_{\mathrm{base}}.
\]
For this numerical illustration only, take $r_A=0.20$ and assume
that $B_n=0.32/\sqrt n$ is a valid RMS reconstruction-displacement bound.
The constant $0.32$ is a hypothetical input, not a quantity inferred from
seismic data or from the preceding elliptic illustration of weight error.
At $n=64$, $B_{64}=0.04$, so the calibrated RMSE upper bound is $0.24$:
this is below a baseline of $0.30$, but does not settle the comparison
with a baseline of $0.22$. The sufficient condition first holds at $n=11$
and $n=257$, respectively, in these hypothetical scenarios
(\cref{fig:seismic-margin}). At $n=256$ the upper bound equals $0.22$,
which does not certify a strict improvement over the narrower baseline.
These calculations show how the oracle margin changes the calibration count
needed to satisfy the sufficient condition. They provide no observed RMSE
improvement, empirical validation of $B_n$, or acquisition recommendation.

To turn this illustration into a wave-inversion test, one would need to
justify a calibration model whose real-symmetric replicates have the
correct mean for the real representation of the wave-data Gram, with
conditional IID sampling, finite second moments, and independence from
the test data conditional on the configuration. Using
\cref{prop:frozen-stability} to obtain $B_n$ would also require a verified
oracle strong-convexity constant on the specified convex search set and
comparison of minimizers on that set. The published reconstructions do
not establish these hypotheses or measure a margin under independent
noisy calibration. This criterion does not ensure that a nonlinear solver
finds those minimizers; the additional acquisition burden of calibration
must also be reported.

\begin{remark}[When does calibration preserve a monotonicity interval?]
\label{rem:calibrated-monotonicity}
Consider the frozen objective forms $F_A$ and
$F_{\widehat A}$ introduced in
\cref{prop:frozen-stability}, but do not impose the
strong-convexity assumption of that proposition.
Work in the scalar-coefficient case, and let $J$ be a
fixed compact interval contained in $\cC_h\subset\R$.
Fix the test data, and use the same residual and
regularization in both objectives. Assume that $e$ and
the common regularizer $\cR$ are $C^1$ on a neighborhood
of $J$. Suppose that
\[
\gamma_J:=\inf_{c\in J}|F_A'(c)|>0,
\qquad
M_J:=\sup_{c\in J}
\norm{e'(c)}_2\,\norm{e(c)}_2<\infty.
\]
Since the weights are frozen,
\[
\sup_{c\in J}|F_{\widehat A}'(c)-F_A'(c)|
\le M_J\norm{\widehat A^{-1}-A^{-1}}_2.
\]
Thus, whenever the right-hand side is strictly smaller
than $\gamma_J$, the calibrated derivative has the same
sign as the oracle derivative throughout $J$. In
particular, calibration preserves strict monotonicity
and introduces no stationary point on $J$.

Under the hypotheses of \cref{cor:calibration-consistency},
set $\widehat K=\widehat K_+$ and $K=K(S)$, and define
\[
A=\Sigma+K(S),\qquad
\widehat A=\Sigma+\widehat K_+.
\]
Take the common regularization to be fixed independently
of the calibration replicates, and
assume that calibration is independent of the test data
conditional on $S,\Sigma$. Then $M_J$ and $\gamma_J$ are
measurable with respect to $\sigma(S,\Sigma,e)$. The derivative
bound and the inequality between the operator and Frobenius
norms give the event inclusion
\[
\left\{\exists\,c\in J:\ F_{\widehat A}'(c)=0\right\}
\subseteq
\left\{M_J^2\norm{\widehat A^{-1}-A^{-1}}_F^2\ge\gamma_J^2\right\}.
\]
By conditional independence, the bound in \cref{eq:inverse-L2}
also holds after conditioning on $e$. Conditional Markov's
inequality therefore gives, on $\{\gamma_J>0\}$,
\[
\mathbb P\!\left(
\exists\,c\in J:\ F_{\widehat A}'(c)=0
\,\middle|\, S,\Sigma,e
\right)
\le
\min\!\left\{1,\frac{M_J^2\,v_n(S)}
{\sigma_0^4\,\gamma_J^2}\right\}.
\]
No division by $M_J$ is needed: on $\{\gamma_J>0\}$,
if $M_J=0$, the event on the right of the inclusion
is empty, and hence the conditional probability is zero.
The result does not apply on an interval containing an
oracle stationary point, where $\gamma_J=0$.
It is neither a global count of local minima nor a convergence
guarantee for an optimization algorithm.
\end{remark}

\paragraph{One-dimensional Helmholtz example.}
We use the constant-speed point-source geometry of Section~3.3.3 of
\cite{vanLeeuwenYang2025}. On $(0,1)$, solve
$-u_j''-(k/c)^2u_j=\delta_{x_j}$ with $u_j(0)=0$ and
$u_j'(1)-\mathrm{i}(k/c)u_j(1)=0$. We take $k=40$, $x_j=j/9$
for $j=1,\dots,8$, and noiseless test data at $c_0=1$.
With $q=k/c$, the solution is
$u_j(x;c)=q^{-1}\sin(q\min(x,x_j))\exp(\mathrm{i}q\max(x,x_j))$.
The columns of $E(c)$ are conjugated receiver responses minus the
corresponding test data. The oracle Gram has entries
$G_{ij}=\int_0^1\overline{u_i'(x;c_0)}u_j'(x;c_0)\,dx$,
evaluated analytically by splitting at the sources. Set $\rho=0.1$.
Use $K=\rho^{-1}\mathcal{T}(G)$, where
$\mathcal{T}(G)=\left(\begin{smallmatrix}\Re G&-\Im G\\\Im G&\Re G\end{smallmatrix}\right)$,
and stack the real and imaginary parts of each residual column in $e(c)$.
The ordinary and oracle objectives are
$\frac12\operatorname{tr}(e(c)^T W e(c))$ with weights $I$ and
$A^{-1}=(I+K)^{-1}$, respectively. For the calibrated objective, write
$\widehat A^{-1}=(I+\widehat K_+)^{-1}$. All use the same residual,
search domain $[0.5,1.5]$, and zero regularization. Both nonidentity
weights are fixed as $c$ varies.

For the synthetic calibration model, $X_s=K+Z_s$ are IID real-symmetric
$16\times16$ matrices. For $n$ independent sessions, define
\[
\bar X_n=\frac1n\sum_{s=1}^n X_s,
\qquad \widehat K_+=[\bar X_n]_+.
\]
Each upper-triangular entry of $Z_s$, including the diagonal, is
independently sampled from $\mathcal N(0,0.12^2)$ and mirrored. The calibration thus has the correct realified mean, finite
second moments, and is independent of the fixed test data. This is noise
on the scaled Gram, not a model of raw sensor data or a reproduction of
the reference paper's frequency-derivative calibration protocol.

In this exploratory configuration, $J=[1.114,1.130]$ was selected by
inspecting the noiseless objectives and then held fixed across all
calibration draws. The ordinary derivative is about $-0.01326$ and
$0.01630$ at the left and right endpoints, so ordinary least squares has
an interior local minimizer on $J$ away from $c_0$. For the frozen oracle,
analytic derivative bounds padded over a grid give
$\gamma_J\ge 0.0622$ and $M_J\le 0.923$ in
\cref{rem:calibrated-monotonicity}, with the Euclidean norm on stacked
residual columns. Specifically, if $B_i$ bounds the norm of the $i$th
residual derivative ($i=0,1,2$), then
$|F_A''|\le\|A^{-1}\|_2(B_2 B_0+B_1^2)$; this controls the derivative
between grid nodes. The computed lower bound for the oracle derivative is positive throughout $J$. These inter-node bounds are obtained from analytically derived derivative estimates evaluated in standard floating-point arithmetic; no validated interval-arithmetic certificate is claimed.

\begin{figure}[tbp]
  \centering
  \includegraphics[width=\linewidth]{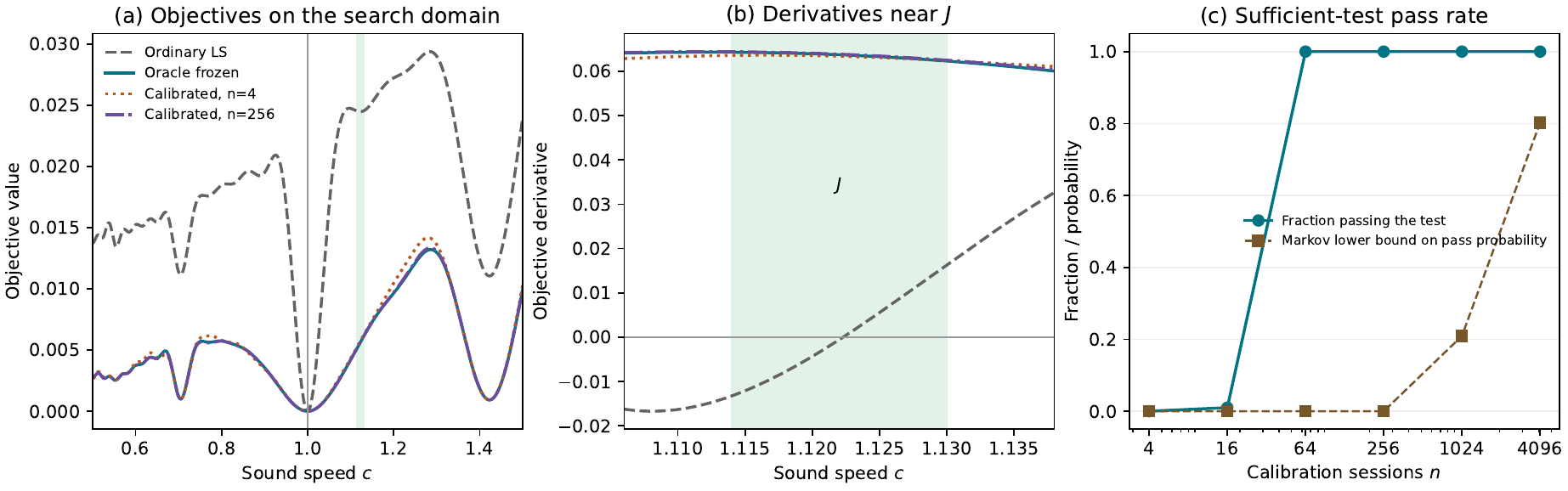}
  \caption{Synthetic 1D Helmholtz experiment. Left: ordinary, oracle-frozen,
  and two calibrated-frozen objectives; the calibrated examples use the first
  draw at each labeled $n$. Center: derivatives near the selected interval
  $J$ (shaded). Right: empirical fraction passing the sufficient test and a
  Markov lower bound on its probability; neither curve reports solver success.}
  \label{fig:helmholtz-calibration-landscape}
\end{figure}

Among $500$ independent calibration draws per session count, the operator-
norm sufficient test passed in $0$, $5$, and $500$ draws for $n=4,16,64$,
respectively, and in all draws at the larger counts tested
(\cref{fig:helmholtz-calibration-landscape}). Failure of the test is
inconclusive: for the first draw at $n=4$, the perturbation criterion does
not pass, but a direct padded derivative calculation still gives a positive
slope on $J$. The same figure displays the conservative Markov lower bound
on passing the sufficient test, using $v_n=16^2(0.12)^2/n$ and $\sigma_0=1$.
The pass counts are not estimates of the probability of a stationary point,
and $500/500$ is not a probability-one result. The conclusions concern this
interval and calibration model only, not global geometry, solver convergence,
reconstruction RMSE, or the cost-effectiveness of additional acquisitions.

\begin{remark}[Preserving the existence of an interior local minimizer under calibration]
\label{rem:local-minimum-persistence}
Use the frozen objective forms in
\cref{prop:frozen-stability}, without its
strong-convexity assumption. Let
$D\subset\operatorname{int}(\cC_h)$ be a fixed compact
set, and let $c_\star\in\operatorname{int}(D)$ be a fixed
comparison point. Assume that $e$ and the common
regularizer $\cR$ are continuous on $D$. Take $D$,
$c_\star$, and the common regularization term
$\varrho\cR$ to be deterministic. Suppose the oracle
objective has a positive boundary gap,
\[
  \delta_D:=\min_{c\in\partial D}
    \bigl(F_A(c)-F_A(c_\star)\bigr)>0.
\]
Define the residual contrast
\[
  L_D:=\frac12\max_{c\in\partial D}
  \norm{e(c)e(c)^T-e(c_\star)e(c_\star)^T}_F.
\]
The common regularization cancels in the perturbation of
each objective difference. By the Frobenius
inner-product inequality,
\[
  \min_{c\in\partial D}
    \bigl(F_{\widehat A}(c)-F_{\widehat A}(c_\star)\bigr)
  \ge \delta_D -
      L_D\norm{\widehat A^{-1}-A^{-1}}_F.
\]
By continuity and compactness, both objectives attain
their minima over $D$. The positive boundary gap implies
that every minimizer of $F_A$ over $D$ lies in
$\operatorname{int}(D)$. If the right-hand side of the
preceding inequality is strictly positive, the same is
true of every minimizer of $F_{\widehat A}$ over $D$.
Each such interior minimizer is a local minimizer of
the corresponding objective on $\cC_h$.

Under the hypotheses of \cref{cor:calibration-consistency},
set $\widehat K=\widehat K_+$ and $K=K(S)$, and define
\[
 A=\Sigma+K(S),\qquad
 \widehat A=\Sigma+\widehat K_+.
\]
Assume that calibration is independent of the test data
conditional on $S,\Sigma$. Then $\delta_D$ and $L_D$
are measurable with respect to $\sigma(S,\Sigma,e)$, and
the bound in \cref{eq:inverse-L2} remains valid after
conditioning on $e$. Conditional Markov's inequality gives
\[
  \mathbb P\!\left(
    \min_{c\in\partial D}F_{\widehat A}(c)
       \le F_{\widehat A}(c_\star)
    \,\middle|\, S,\Sigma,e
  \right)
  \le \min\left\{1,
    \frac{L_D^2 v_n(S)}
         {\sigma_0^4\delta_D^2}
  \right\}
  \qquad\text{on }\{\delta_D>0\}.
\]
Indeed, failure of the strict boundary comparison implies
\[
 L_D^2\norm{\widehat A^{-1}-A^{-1}}_F^2
 \geq \delta_D^2.
\]
If $L_D=0$ and $\delta_D>0$, failure is impossible.
The bound concerns failure of a sufficient condition,
not uniqueness, the location of a particular minimizer,
or convergence of a solver.
\end{remark}

\subsection{Finite-sample bias and remaining statistical questions}\label{sec:calibration-scope}
Projection can introduce bias and inversion is nonlinear. Even an unbiased positive-definite estimate need not give an unbiased inverse-weighted quadratic. \Cref{prop:random-weight,cor:inverse-bias} state the relevant conditional identity and a small-perturbation expansion. These finite-sample facts are fully compatible with \cref{cor:calibration-consistency}; bias is not the same question as consistency.

The theorem requires IID sessions and a correct matrix mean. Dependent or drifting acquisitions, excitation uncertainty, misspecified symmetrization, and fitting $S$ on the test calibration sample are not covered. Independent non-identically distributed replicates with a common mean still give unbiasedness, but require index-dependent variance terms rather than \cref{eq:calibration-variance}. No claim of efficient calibration design or superiority of the PSD estimator over shrinkage or structured alternatives is made.

\section{Limitations and open questions}\label{sec:open}
\paragraph{Acquisition and geometry.}
The reciprocal identity is proved for bounded point observation in a one-dimensional, symmetric, uniformly elliptic state model. Extension to higher dimensions or non-reciprocal operators requires different observation/source spaces and a new protocol. The fixed-reference construction keeps its changing geometry explicit; it does not validate Green identities for unrelated forcing norms. Calibration at $c_0$ is not transport to another coefficient. A function-valued calibration $c\mapsto K(c)$ needs additional data and modeling.

\paragraph{Reconstruction and data budgets.}
The included numerical illustrations use synthetic calibration noise or hypothetical RMSE inputs; they are not reconstruction experiments and provide no evidence that a norm removes nonlinear local minima, enlarges convergence basins, or improves efficiency. \Cref{prop:frozen-stability} assumes coefficient-level strong convexity and controls only distance to the oracle frozen estimator. Future experiments must include unfavorable regimes, fair optimization baselines, and explicit accounting for the matrix calibration acquisitions and any extra frequencies or training data.

\paragraph{Learning.}
Compactness proves finite-dimensional optimistic attainment even with nonunique lower solutions; it does not prove population consistency, uniqueness of the upper optimum, or convergence of practical algorithms. Local hypergradients do not cover branch switching, nonsmooth penalties without a reformulation, spectral zero crossings, or transitions between competing minima. Function-space bilevel limits require separate compactness and lower-semicontinuity arguments.

\paragraph{Probability.}
A learned deterministic geometry need not be a physical covariance. When both measurement and source covariances are unrestricted, their sum generally does not identify the two components separately. The trace-class distinction and range-domain factor formula do not amount to a full infinite-dimensional Bayesian analysis of the nonlinear coefficient problem. Concentration under dependence, calibrated confidence regions, and distribution shift remain open here. In contrast, fixed-target mean-square calibration and inverse-weight consistency are already established by \cref{cor:calibration-consistency}.

\section{Conclusion}\label{sec:conclusion}
Quadratic source elimination exposes a measurement-space compression, but its deterministic and probabilistic meanings must be stated separately. For the structured target $G+GSG$, independent cross-replicates yield an explicit unbiased estimator and exact finite-sample variance. PSD projection supplies an admissible frozen weight, and the resulting mean-square bound transfers to inverse weights and, under strong convexity, to the oracle frozen coefficient estimator.

The reciprocal example supplies a concrete but restricted measurement protocol and shows why the source/observation distinction matters. Exact candidate-dependent relaxation, frozen calibrated reconstruction, and supervised training are distinct operations. Their separation makes the claims auditable without asserting a new general elimination principle, physical covariance recovery, or universal nonlinear reconstruction improvement.

\appendix
\section{Hilbert-valued Gaussian covariance}\label{app:geometry}
\begin{proposition}[Coercivity versus trace class]\label{prop:trace-class}
Let $\cQ$ be separable and infinite-dimensional. A bounded self-adjoint operator $C\succeq\kappa I$ with $\kappa>0$ is not trace class. A bounded positive self-adjoint operator is the covariance of a centered $\cQ$-valued Gaussian random element if and only if it is trace class; then $\E\norm{q}^2=\tr C$.
\end{proposition}
\begin{proof}
For any orthonormal basis $(\varphi_k)$,
\[
 \sum_{k=1}^N\ip{C\varphi_k}{\varphi_k}\ge N\kappa,
\]
so a coercive positive operator has infinite trace. For necessity of trace class in the Gaussian assertion, a Hilbert-valued Gaussian has finite second moment by the standard Gaussian integrability theorem (Fernique; see \cite{Bogachev1998,DaPratoZabczyk2014}). Parseval and Tonelli in an \emph{arbitrary} orthonormal basis give
\[
 \sum_k\ip{C\varphi_k}{\varphi_k}
 =\sum_k\E\abs{\ip{q}{\varphi_k}}^2=\E\norm{q}^2<\infty.
\]
No eigenbasis of a general bounded positive operator is assumed.

Conversely, positive trace class implies compactness and permits an orthonormal eigenbasis, including the kernel. If $C\varphi_k=\rho_k\varphi_k$ and $\xi_k$ are independent standard normals, then
\[
 q=\sum_k\sqrt{\rho_k}\,\xi_k\varphi_k
\]
converges in mean square because $\sum_k\rho_k<\infty$. It also converges almost surely: the expected sum of the squared coefficients is finite, so that sum is finite almost surely. Every continuous linear functional of the limit is Gaussian by convergence of the corresponding scalar Gaussian series. The covariance is $C$, and $\E\norm{q}^2=\sum_k\rho_k=\tr C$.
\end{proof}

\section{Local statistical identities}\label{app:statistics}
These are classical finite-dimensional reference identities, not nonlinear PDE recovery results. Let $J\in\R^{m\times p}$ and let $\Gamma^\dagger\succ0$ denote a true local noise covariance. The superscript $\dagger$ in this appendix is a label for the true covariance, not a generalized inverse.

\subsection{Gaussian posterior and risk}\label{sec:statistics-posterior}
\begin{proposition}[Gaussian posterior mean and exact risk]\label{prop:posterior}
Let $h\sim\mathcal N(0,P)$ and $\zeta\sim\mathcal N(0,\Gamma^\dagger)$ independently, with $P\succ0$, and let $y=Jh+\zeta$. Define
\[
 V=(P^{-1}+J^T(\Gamma^\dagger)^{-1}J)^{-1},\qquad
 m(y)=VJ^T(\Gamma^\dagger)^{-1}y=PJ^T(JPJ^T+\Gamma^\dagger)^{-1}y.
\]
Then $h\mid y\sim\mathcal N(m(y),V)$, and for every square-integrable estimator $g(y)$,
\begin{equation}\label{eq:posterior-risk}
 \E\norm{h-g(y)}_2^2=\tr V+\E\norm{m(y)-g(y)}_2^2.
\end{equation}
Thus the posterior mean minimizes squared-error risk.
\end{proposition}
\begin{proof}
Complete the square in $h$ in the joint Gaussian density. For the second form of $m$, let $H=P^{-1}+J^T(\Gamma^\dagger)^{-1}J$ and $B=JPJ^T+\Gamma^\dagger$. The identity $HPJ^T=J^T(\Gamma^\dagger)^{-1}B$, multiplied on the left by $H^{-1}$ and the right by $B^{-1}$, gives the result. Finally $\E[h-m(y)\mid y]=0$ makes the cross term in $h-g=(h-m)+(m-g)$ vanish, giving \cref{eq:posterior-risk}.
\end{proof}
The quadratic convention is $\tfrac12\norm{y-Jh}_{(\Gamma^\dagger)^{-1}}^2+\tfrac12h^T\Lambda h$ with $\Lambda=P^{-1}$. A penalty $\lambda\norm{h}^2$ corresponds to precision $2\lambda I$. If the data term is multiplied by a positive scalar $a$, the penalty precision must also be multiplied by $a$ to preserve the estimator.

\subsection{Regularized linear estimators and genuine LMMSE}\label{sec:statistics-lmmse}
\begin{proposition}[Linear error covariance and LMMSE]\label{prop:lmmse}
Let $h,\zeta$ be centered, square-integrable, and uncorrelated, with $\Cov(h)=P\succeq0$, $\Cov(\zeta)=\Gamma^\dagger\succ0$, and $y=Jh+\zeta$. For $\Gamma\succ0$, $\Lambda\succeq0$, define
\[
 M_\Gamma=J^T\Gamma^{-1}J,\quad N_\Gamma=J^T\Gamma^{-1}\Gamma^\dagger\Gamma^{-1}J,\quad
 H_\Lambda=M_\Gamma+\Lambda.
\]
Assume $H_\Lambda\succ0$. The regularized linear estimator
\[
 \widehat h_\Lambda=H_\Lambda^{-1}J^T\Gamma^{-1}y
\]
has centered unconditional error with covariance
\begin{equation}\label{eq:ridge-cov}
 \Cov(\widehat h_\Lambda-h)=H_\Lambda^{-1}(N_\Gamma+\Lambda P\Lambda)H_\Lambda^{-1}.
\end{equation}
For arbitrary weights it need not be LMMSE. The unique squared-risk optimal linear map is
\begin{equation}\label{eq:true-lmmse}
 A_\star=PJ^T(JPJ^T+\Gamma^\dagger)^{-1}.
\end{equation}
If $P\succ0$, the matched choices $\Gamma=\Gamma^\dagger$, $\Lambda=P^{-1}$ give $\widehat h_\Lambda=A_\star y$. Zero unconditional mean error is not fixed-parameter or conditional unbiasedness.
\end{proposition}
\begin{proof}
The ridge error equals $H_\Lambda^{-1}(J^T\Gamma^{-1}\zeta-\Lambda h)$. Uncorrelatedness gives \cref{eq:ridge-cov}. For $B=\Cov(y)=JPJ^T+\Gamma^\dagger\succ0$, $A_\star B=PJ^T=\Cov(h,y)$, so
\[
 \E\norm{Ay-h}^2=\tr(P-PJ^TB^{-1}JP)
                  +\norm{(A-A_\star)B^{1/2}}_F^2.
\]
This proves unique linear optimality without Gaussianity \cite{Kay1993}. The matched normal-equation form follows from the same push-through identity as in \cref{prop:posterior}.
\end{proof}
The specialization $\Lambda=0$ is defined only when $J$ has full column rank. When $\Gamma=\Gamma^\dagger$ it yields covariance $M_\Gamma^{-1}$. Its deterministic-parameter optimality is the separate BLUE statement below, not a claim that unregularized GLS is LMMSE for a random coefficient prior.

\subsection{GLS/BLUE for a deterministic local coefficient}\label{sec:statistics-gls}
\begin{proposition}[Generalized least squares]\label{prop:gls}
Let $h\in\R^p$ be deterministic, $J$ have full column rank, and $y=Jh+\xi$ with centered square-integrable $\xi$ of covariance $\Gamma^\dagger\succ0$. Then
\[
 \widehat h_{\mathrm{GLS}}=(J^T(\Gamma^\dagger)^{-1}J)^{-1}J^T(\Gamma^\dagger)^{-1}y
\]
is unbiased with covariance $(J^T(\Gamma^\dagger)^{-1}J)^{-1}$ and has minimal covariance in Loewner order among linear unbiased estimators.
\end{proposition}
\begin{proof}
Let its matrix be $A_0$, so $A_0J=I$. Any other unbiased linear map is $A_0+B$ with $BJ=0$. The covariance difference is $B\Gamma^\dagger B^T\succeq0$, since $A_0\Gamma^\dagger B^T=(J^T(\Gamma^\dagger)^{-1}J)^{-1}J^TB^T=0$.
\end{proof}
For independent or uncorrelated centered source and measurement errors $\xi=\varepsilon+Tq$, the true covariance is $\Gamma^\dagger=\Sigma+TC_{\mathrm{prob}}T^*$. Equality of this physical covariance with a learned deterministic weight is an additional assumption, not a consequence of training.

\subsection{Random inverse-weighted losses}\label{sec:statistics-random-weight}
\begin{proposition}[Conditional quadratic expectation]\label{prop:random-weight}
Let $Y\in\R^m$ be square-integrable and let $\widehat A$ be a measurable random SPD matrix, with $\E[Y^T\widehat A^{-1}Y]<\infty$. Then almost surely,
\[
 \E[Y^T\widehat A^{-1}Y\mid\widehat A]
 =\tr\{\widehat A^{-1}\Cov(Y\mid\widehat A)\}
       +\E[Y\mid\widehat A]^T\widehat A^{-1}\E[Y\mid\widehat A].
\]
\end{proposition}
\begin{proof}
On the events $\{\norm{\widehat A^{-1}}_2\le N\}$ the conditional expansion of the square is integrable, and the centered cross term has conditional expectation zero. These events exhaust the sample space, giving the almost-sure identity by localization. The assumed finite weighted moment ensures integrability of the nonnegative terms on the right.
\end{proof}
\begin{corollary}[Small-perturbation inverse bias]\label{cor:inverse-bias}
Let $A\succ0$, let $\Delta$ be symmetric with $\E\Delta=0$ and $\norm{\Delta}_2\le\lambda_{\min}(A)/2$ almost surely, and let $e$ be deterministic. Set $v=A^{-1/2}e$ and $B=A^{-1/2}\Delta A^{-1/2}$. Then
\[
 \E[e^T(A+\Delta)^{-1}e]=e^TA^{-1}e+\E\norm{Bv}^2+R,
 \qquad |R|\le2\norm{v}^2\norm{A^{-1}}_2^3\E\norm{\Delta}_F^3.
\]
In particular,
\begin{equation}\label{eq:inverse-bias-bound}
 0\le\E[e^T(A+\Delta)^{-1}e]-e^TA^{-1}e
 \le2\norm{v}^2\norm{A^{-1}}_2^2\E\norm{\Delta}_F^2.
\end{equation}
\end{corollary}
\begin{proof}
Use $(I+B)^{-1}=I-B+B^2-B^3(I+B)^{-1}$ and $\norm{B}_2\le1/2$. The linear term vanishes in expectation, and the remainder bound follows from $\norm{(I+B)^{-1}}_2\le2$. For \cref{eq:inverse-bias-bound}, the exact identity $(I+B)^{-1}=I-B+B^2(I+B)^{-1}$ shows that the bias is $\E[v^TB^2(I+B)^{-1}v]$, which lies between $0$ and $2\E\norm{Bv}^2$.
\end{proof}
This expansion assumes an unbiased, uniformly small perturbation. It is not automatically an expansion for the PSD-projected estimator, whose finite-sample mean can differ from its target. The separate consistency bound \cref{eq:inverse-L2} does apply to that estimator.

\section{Proof of the calibration theorem}\label{app:calibration}
Condition throughout on an independently fixed $S$. The matrices $X_i$ remain IID symmetric with mean $G$ and finite second moments. Write $W_i=X_i-G$. Independence gives $\E[X_iSX_j\mid S]=GSG$ for $i\ne j$, proving unbiasedness. A useful implementation identity is
\begin{equation}\label{eq:calibration-scatter}
 \widehat K_U=\bar X+\bar XS\bar X-\frac1n\widehat B_S,\qquad
 \widehat B_S=\frac1{n-1}\sum_i(X_i-\bar X)S(X_i-\bar X).
\end{equation}
It follows by expanding the ordered-pair sum; no commutativity is needed.

Expanding around the true mean gives the Hoeffding decomposition
\begin{equation}\label{eq:calibration-decomp}
 \widehat K_U-K(S)=\frac1n\sum_iL_i+\frac2{n(n-1)}\sum_{i<j}H_{ij},
\end{equation}
where $L_i=\mathcal L_S(W_i)$ and $H_{ij}=\mathcal H_S(W_i,W_j)$. The linear terms are centered. Each pair kernel is conditionally degenerate:
\[
 \E[H_{ij}\mid W_i,S]=\E[H_{ij}\mid W_j,S]=0.
\]
The bounds
\[
 \norm{\mathcal L_S(W)}_F\le(1+2\norm{G}_2\norm{S}_2)\norm{W}_F,\qquad
 \norm{\mathcal H_S(W,W')}_F\le\norm{S}_2\norm{W}_F\norm{W'}_F
\]
show square integrability. In the second bound, independence factors the squared product moment, so a fourth moment of one replicate is not needed.

Distinct linear terms are orthogonal in $L^2$. Linear--pair terms vanish by independence when disjoint and by conditional degeneracy when sharing an index. Disjoint pair terms are independent and centered. For two pairs sharing exactly one index, condition on the shared error and on $S$; the remaining errors are conditionally independent and each kernel is conditionally centered, so their Frobenius inner product has zero conditional expectation. \emph{Conditional independence as well as centering} is used in this last step.

Consequently the variance of \cref{eq:calibration-decomp} is the sum of $n$ equal linear second moments and $\binom n2$ equal pair second moments:
\[
 \Var_F(\widehat K_U\mid S)
 =\frac{n}{n^2}a(S)
  +\frac{4}{n^2(n-1)^2}\binom n2 b(S)
 =\frac{a(S)}n+\frac{2b(S)}{n(n-1)}.
\]
This proves \cref{thm:calibration}. IID is used for the common constants, independence for orthogonality, and sample independence of $S$ for the conditional replicate law. These hypotheses are not replaced by symmetrizing dependent acquisitions.

\section{Supplementary reproducibility checks}\label{app:verification}
The accompanying verification files \texttt{verify\_revised\_v4.py} and \texttt{legacy\_algebra\_checks.py} contain the numerical material. The latter retains the earlier representer/Monte Carlo checks; the former includes deterministic finite-support and derivative checks. The arXiv ancillary verification run used Python~3.13.5 and NumPy~2.3.5. See \texttt{README\_v4.md} and \texttt{requirements\_v4.txt} for setup. From the directory containing the verification files, run
\begin{verbatim}
OPENBLAS_NUM_THREADS=1 python verify_revised_v4.py
\end{verbatim}
The saved output is \texttt{verification\_revised\_v4.txt}. In particular:
\begin{itemize}[leftmargin=*]
\item A noncommuting $2\times2$ finite-support example enumerates all $64$ outcomes of three IID replicates. It gives exact/predicted Frobenius MSE $1.557914666666667$, up to floating-point rounding, and mean error below $10^{-12}$.
\item The IID scalar example in \cref{rem:counterexample} is fully enumerated, checking both $\E\widehat K_U=2$ and $\E\widehat K_+=3$.
\item Fixed-reference finite-dimensional elliptic analogues check $MC_0M^*=M\cA^{-1}M'$, the state-geometry compression, and the full coefficient derivative by finite differences. The forcing-map counterexample is retained as a model-distinction check, not as a failure of the corrected lemma.
\item Singular-factor uniqueness, the whitened stability constant, weighted push-through identities, LMMSE optimality, calibration projection/inverse bounds, and the frozen coefficient perturbation bound are checked on explicit finite-dimensional instances.
\end{itemize}
These are audit aids, not replacements for proofs or evidence of reconstruction performance. Last digits and Monte Carlo output can vary with NumPy/BLAS versions. No empirical PDE reconstruction claim is made.
The three numerical illustrations are accompanied by their vector figures in the manuscript source package and by source scripts, generated data, runtime requirements, and verification notes as ancillary material. The elliptic and Helmholtz scripts simulate synthetic matrix-level calibration; the seismic-margin script evaluates a hypothetical sufficient-condition calculation rather than simulating an inversion.

\end{document}